\documentclass[11pt,twoside]{article}
\usepackage{times}
\usepackage{amsmath,amssymb}
\usepackage{color}
\usepackage{graphicx}

\allowdisplaybreaks
 \def \no{\nonumber}

\newcommand{\be}{\begin{equation}}
\newcommand{\ee}{\end{equation}}
\newcommand{\bea}{\begin{eqnarray}}
\newcommand{\eea}{\end{eqnarray}}

\def\R {\Bbb R}

\def\p{\partial}
\def\ve{\varepsilon}

\def\la{\lambda}

\def\al{\alpha}

\def\t{\tilde}

\def\vp{\varphi}
\def\O{\Omega}

\def\p{\partial}

\def\ve{\varepsilon}

\def\la{\lambda}
\def\al{\alpha}

\def\t{\tilde}

\def\vp{\varphi}

\def\ka{\kappa}

\def\ds{\displaystyle}

\def\t{\tilde}
\def\vp{\varphi}

\def\ga{\gamma}

\def\e{\epsilon}

\def\be{\begin{equation}}
\def\ee{\end{equation}}
\def\bes{\begin{equation*}}
\def\ees{\end{equation*}}
\def\ba{\begin{aligned}}
\def\ea{\end{aligned}}
\def\wmp{w^{m+1}}
\def\wm{w^m}
\def\wmm{w^{m-1}}
\def\gm{g^{m}}
\def\gmp{g^{m+1}}

\def\bw{\bar{w}}
\def\br{\bar{\rho}}
\def\bv{\bar{v}}
\def\bt{\bar{\theta}}

\def\Lat{\Lambda_{\tau}}
\def\pb{\phi_{\beta}}

\def\pt{\p_{\tau}}

\def\sm{\setminus}
\def\mut{\tilde{\mu}}
\def\bal{\bar{\alpha}}
\def\tal{\tilde{\alpha}}
\def\hal{\hat{\alpha}}
\def\kp{k_{\phi}}
\def\les{\lesssim}
\def\Kp{K_{\phi}}
\def\Gp{\Gamma_{\phi}}

\def\AehN{A_{\ka,h,N}}

\def\AehhN{A_{h^2\ka/N,h,N}}
\def\de{\delta}

\def\tn{\tau_{n}}
\def\tnp{\tau_{n+1}}
\def\msV{\mathscr{V}}
\def\msL{\mathscr{L}}

\def\vphi{\varphi}

\def\iOR{\int_{\Omega\times\R^3}}
\def\ep{\varepsilon}
\def\iO{\int_{\O}}
\def\tmn{\tilde{\mu}^{-1/2}}
\def\tmp{\tilde{\mu}^{1/2}}
\def\Lap{\Delta}
\def\epn{\frac{1}{\varepsilon}}
\def\ist{\int_{s}^{\tau}}

\def\msV{\mathcal{V}}
\def\msL{\mathcal{L}}
\allowdisplaybreaks

\begin{document}
 \footskip=0pt
 \footnotesep=2pt
\let\oldsection\section
\renewcommand\section{\setcounter{equation}{0}\oldsection}
\renewcommand\thesection{\arabic{section}}
\renewcommand\theequation{\thesection.\arabic{equation}}
\newtheorem{claim}{\noindent Claim}[section]
\newtheorem{theorem}{\noindent Theorem}[section]
\newtheorem{lemma}{\noindent Lemma}[section]
\newtheorem{proposition}{\noindent Proposition}[section]
\newtheorem{definition}{\noindent Definition}[section]
\newtheorem{remark}{\noindent Remark}[section]
\newtheorem{corollary}{\noindent Corollary}[section]
\newtheorem{example}{\noindent Example}[section]

\title{The global existence and large time behavior of smooth compressible fluid  in an infinitely  expanding ball, III:
\protect\\ the 3-D Boltzmann equation}

\author{\quad
Huicheng Yin$^{1,*}$, \quad Wenbin Zhao$^{2}$\footnote{*Huicheng Yin (huicheng$@$nju.edu.cn, 05407@njnu.edu.cn)
and Wenbin Zhao (zhaowb1989@gmail.com) are
supported by the NSFC (No.11571177) and A Project Funded by the Priority Academic Program Development of
Jiangsu Higher Education Institutions. }\vspace{0.5cm}\\
\small 1. School of Mathematical Sciences, Jiangsu Provincial Key Laboratory for Numerical Simulation\\
\small of Large Scale Complex Systems, Nanjing Normal University, Nanjing 210023, China.\\
\small 2. Department of Mathematics and IMS, Nanjing University, Nanjing 210093, China\\
}

\date{}
\maketitle
% \vskip 0.2in

\centerline {\bf Abstract} \vskip 0.3 true cm

This paper is a continuation of the works in \cite{Euler} and \cite{NS}, where the authors have established
the global existence of smooth compressible flows in infinitely expanding balls
for inviscid gases and viscid gases, respectively.
In this paper, we are concerned with the global existence and large time behavior
of compressible Boltzmann gases in an infinitely expanding ball. Such a problem is one of
the interesting  models in studying the theory of global smooth solutions to multidimensional
compressible gases with time dependent boundaries and vacuum states at infinite time.
Due to  the  conservation of mass, the fluid in the  expanding ball becomes rarefied and eventually
tends to a vacuum state meanwhile  there
are no appearances of vacuum domains in any part of the expansive ball,
which is easily observed in finite time. In the present paper,
we will confirm this physical phenomenon for the Boltzmann equation
by  obtaining the exact lower and upper bound on the macroscopic density function.

\vskip 0.3 true cm

\vskip 0.3 true cm


\section{Introduction}

The compressibility of gases plays a basic role in gas dynamics. When
one squeezes a soft container filling with gases, the gases will  become denser and the
corresponding temperature will get higher in the adiabatic process.
In this paper, as in \cite{Euler} and \cite{NS}, we consider an opposite situation for the compressible gases
filling a 3-D expansive ball. It is assumed that the expansive ball
is described by
\bes
  \Omega_t=\{x=(x_1,x_2,x_3)\in\R^3:\, |x|=\sqrt{x_1^2+x_2^2+x_3^2}< R(t), t\ge 0\},
\ees
at
the time $t$, where $R(t)=(1+h^2 t^2)^{1/2}$ for some positive constant $h$. From the expression of $\O_t$,
we know that the expansive ball at time $t$ is
formed by pulling out the initial unit ball $\Omega_0=\{x: |x|<1\}$ with smooth speed and acceleration
(see Figure 1 below).
The pulling speed on the boundary is $R'(t)=h^2t(1+h^2t^2)^{-1/2}$, which increases smoothly from $0$ to $h$.
We denote the time-space domain by $S=\{(t,x):\,t>0, |x|<R(t)\}$. Suppose that the movement of the gases in $\Omega_t$ is described by the 3-D Boltzmann equation:
\begin{equation}\label{Boltzmann eq}
\p_t f+\xi\cdot\nabla_x f=Q(f,f),
\end{equation}
where $f=f(t,x,\xi)$ stands for the distribution function of gas particles at time $t$, position $x\in \O_t$ and
velocity $\xi\in\R^3$, the collision operator $Q(f,g)$ with hard-sphere interaction is given by
\begin{equation}\label{Q(f,g)}
Q(f,g)=Q(f,g)(t,x,\xi)=\frac{1}{2}\int_{\R^3\times S^2}
|(\xi-\xi_*)\cdot\omega|(f'g'_*+f'_*g'-fg_*-f_*g)d\xi_*d\omega
\end{equation}
with $\omega\in S^2$ being the unit sphere in $\R^3$, and
\begin{equation*}
    f_*=f(t,x,\xi_*),\,\,\,\,
    f'=f(t,x,\xi'),\,\,\,\,f'_*=f(t,x,\xi'_*),
\end{equation*}
\begin{equation*}
    \xi'=\xi-[(\xi-\xi_*)\cdot\omega]\omega,\,\,\,\,
    \xi'_*=\xi_*+[(\xi-\xi_*)\cdot\omega]\omega.
\end{equation*}

\begin{figure}[htbp]
\centering\includegraphics[width=9cm,height=7cm]{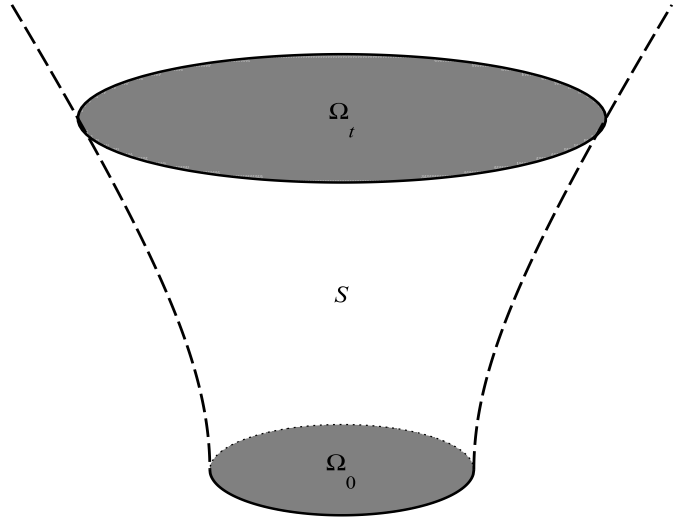}
\caption{\bf A rarefied gas flow in an expanding ball}\label{fig:1}
\end{figure}

\vskip 0.2 true cm
In view of the physical property for the gas flow in $\Omega_t$,
it is plausible to pose the following initial-boundary conditions for equation (\ref{Boltzmann eq}),
\begin{equation}\label{initial-boundary condition}
\left\{ \enspace
\begin{aligned}
f(0,x,\xi) & =f_0(x,\xi), &  &\text{$x\in\O_0$, $\xi\in\R^3$},\\
f(t,x,\xi) & =f(t,x,\xi-2(\xi\cdot n_x)n_x+2R'(t)n_x),
&  & \text{$x\in\p \O_t$, $\xi\cdot n_x<R'(t)$},
\end{aligned}
\right.
\end{equation}
where $n_x$ is the unit outer normal direction of $\{x\in\Bbb R^3: |x|=R(t)\}$ in $\Bbb R^3$,
$R'(t)n_x$ is the velocity of the expansive boundary. Note that the boundary value condition
in \eqref{initial-boundary condition} just only corresponds to the specular-reflection boundary condition.
It is easy to check that the ``traveling Maxwellian''
\begin{equation}\label{traveling Maxwellian}
M(t,x,\xi)=(2\pi)^{-3/2}
\exp\left(-\frac{|\xi|^2}{2}-\frac{h^2}{2}|x-t\xi|^2\right)
\end{equation}
is a special solution to equation (1.1) despite the initial data.

To solve problem \eqref{Boltzmann eq} together with \eqref{initial-boundary condition}, we first make a change of variables $(t,x,\xi)$
such that the expansive domain
$\Omega_t$ becomes a fixed domain $\O=\{y=(y_1, y_2,y_3)\in\Bbb R^3: |y|<1\}$. For this purpose, we set
\begin{equation}\label{A-0}
\left\{ \enspace
\begin{aligned}
\tau & =\frac{1}{h}\arctan(ht), \\
y& =\frac{x}{R(t)}, \\
\eta & =R(t)\xi-\frac{h^2t}{R(t)}x.
\end{aligned}
\right.
\end{equation}
In this case,  $t\in(0,\infty)$ turns into a finite interval $\tau\in(0,\pi/2h)$,
and the special solution $M=M(t,x,\xi)$ becomes
\begin{equation}\label{1-0}
M=M(\tau,y,\eta)=(2\pi)^{-3/2}
\exp\left(-\frac{|\eta|^2+h^2|y|^2}{2}\right).
\end{equation}

Denote by $\mu=\mu(\eta):=(2\pi)^{-3/2}e^{-|\eta|^2/2}$ and $\t\mu=\t\mu(y):=e^{-h^2|y|^2/2}$, we then have $M=\mu\t\mu$
for \eqref{1-0}.
Under the new coordinates $(\tau, y, \eta)$, Boltzmann equation (\ref{Boltzmann eq}) becomes
\begin{equation}\label{eq: f}
\p_{\tau} f+\eta\cdot\nabla_{y}f-h^2y\cdot\nabla_{\eta}f=\cos^2(h\tau) Q(f,f),
\end{equation}
and the initial-boundary data \eqref{initial-boundary condition} become
\begin{equation}\label{data: f}
\left\{ \enspace
\begin{aligned}
f(0,y,\eta)&=f_0(y,\eta), &
& \text{$y\in \O$, $\eta\in\R^3$},\\
f(\tau,y,\eta)&=f(\tau,y,\eta-2(\eta\cdot n_y)n_y), &
& \text{$y\in\p\O$, $\eta\cdot n_y<0$}.
\end{aligned}
\right.
\end{equation}
Note that \eqref{eq: f} is a kind of Boltzmann equation containing a potential term. Denote the transport operator by
\begin{equation}\label{A-A}
\Lambda_{\tau}:=\p_{\tau}+\eta\cdot\nabla_{y}-h^2y\cdot\nabla_{\eta}.
\end{equation}
As in \cite{MR882376}, we define the standard perturbation around $M$ of $f$ by
\begin{equation*}
f=M+M^{1/2}u.
\end{equation*}
In this case, one obtains the Boltzmann equation of $u$ as follows:
\begin{equation}\label{eq for u}
\Lambda_{\tau}u=\t\mu \cos^2(h\tau)Lu+{\t\mu}^{1/2} \cos^2(h\tau)\Gamma(u,u),
\end{equation}
where $L$ and $\Gamma$ are the Boltzmann operators in a bounded domain that can be
expressed as
\begin{equation}\label{L}
Lu=2\mu^{-1/2}Q(\mu,\mu^{1/2}u),
\end{equation}
\begin{equation}\label{Gamma}
\Gamma(u,u)=\mu^{-1/2}Q(\mu^{1/2}u,\mu^{1/2}u).
\end{equation}
Correspondingly, the initial-boundary data of $u$ are
\begin{equation}\label{data for u}
\left\{ \enspace
\begin{aligned}
u(0,y,\eta)&=u_0(y,\eta), &
& \text{$y\in \O$, \quad $\eta\in\R^3$}, \\
u(\tau,y,\eta)&=u(\tau,y,\eta-2(\eta\cdot n_y)n_y), &
& \text{$y\in\p\O$, \quad $\eta\cdot n_y<0$}.
\end{aligned}
\right.
\end{equation}
Suppose that the initial perturbation satisfies
the following conservations:
\be\label{conservation of initial mass}
\int_{\Omega\times\R^3}u_0(y,\eta)M^{1/2}dyd\eta=0,
\ee
\be\label{conservation of initial energy}
\int_{\Omega\times\R^3}(|\eta|^2+h^2|y|^2)u_0(y,\eta)M^{1/2}dyd\eta=0,
\ee
and
\be\label{conservation of initial angular}
\int_{\Omega\times\R^3}(y\times\eta)u_0(y,\eta)M^{1/2}dyd\eta=0.
\ee
Since the specular-reflection boundary conserves both mass, energy and angular momentum,
as in \cite{MR2679358}, without loss of generality, we may assume that the mass, energy
and angular momentum conservation laws hold for all the time. That is, for all $\tau\in[0,\pi/2h)$,
\be\label{law of mass}
\int_{\O\times\R^3}u(\tau,y,\eta)M^{1/2}dyd\eta=0,
\ee
\be\label{law of energy}
\int_{\O\times\R^3}(|\eta|^2+h^2|y|^2)
u(\tau,y,\eta)M^{1/2}dyd\eta=0,
\ee
\be\label{law of angular}
\int_{\O\times\R^3}(y\times\eta)
u(\tau,y,\eta)M^{1/2}dyd\eta=0.
\ee
In order to state our results conveniently, we introduce the following weight function for $\beta>3/2$,
\begin{equation*}
\pb(y,\eta)=(1+|\eta|^2+h^2|y|^2)^{\beta/2}.
\end{equation*}
The main theorem in this paper is

{\bf Theorem 1.1.} {\it For  small $\epsilon_0>0$, suppose that the initial data $f_0=M+M^{1/2}u_0\geq0$ satisfying
(\ref{conservation of initial mass})-(\ref{conservation of initial angular}) with
$\|\pb u_0\|_{\infty}\leq\epsilon_0$, then there exists a constant $\la>0$ and a unique mild solution
$f=M+M^{1/2}u\geq0$ to problem (\ref{eq for u}) together with (\ref{data for u}) and
(\ref{law of mass})-(\ref{law of angular}) such that for $\tau\in(0,\pi/2h)$,
\begin{equation*}
\|\pb u(\tau)\|_{\infty}\leq C e^{-\la\tau}\|\pb u_0\|_{\infty},
\end{equation*}
where $C>0$ is a constant independent of $\tau$. Moreover, if the initial data $f_0(y,\eta)$ is continuous
except on $\ga_{01}$, then $f(\tau,y,\eta)$ is continuous in $[0,\pi/2h)\times\{\bar{\O}\times\R^3\setminus\ga_{01}\}$, where $\gamma_{01}=\{(y,\eta)\in\ga_0: |\eta|\geq h\}$
with $\ga_0=\{(y,\eta)\in\p\Omega\times\R^3:\eta\cdot n_y=0\}$ and $\p\Omega=\{x\in\Bbb R^3: |x|=1\}$.
}

\vskip 0.2 true cm

{\bf Remark 1.1.} {\it
We can also consider the case  that the expanding
speed of the ball is exactly the constant number $h>0$ in Theorem 1.1.
At this time, the radius of the expanding ball at time $t$ is $R_1(t)=1+ht$. Correspondingly,
$M_1(t,x,\xi)=e^{-|\xi-h(x-t\xi)|^2}$ is the background solution of problem (\ref{eq for u}). As long as
we modify the pulling speed near the time $t=0$ to let the speed
increase smoothly from $0$ to $h$, then we can obtain the analogous result to
Theorem 1.1.}

\vskip 0.2 true cm

Next, we study the global physical phenomenon of problem (\ref{eq: f}) together with
(\ref{data: f}). Return to the original coordinates $(t,x,\xi)$ and equation \eqref{Boltzmann eq}.
Let $\rho=\rho(t,x)$, $v=v(t,x)$
and $\theta=\theta(t,x)$ be the mass density, velocity and temperature of the gases, respectively,  i.e,
\begin{equation*}
\begin{split}
&\rho(t,x)=\int_{\R^3}f(t,x,\xi)d\xi,\\
&\left(\rho v\right)(t,x)=\int_{\R^3}\xi f(t,x,\xi)d\xi,\\
&\left(\rho\theta\right)(t,x)=\int_{\R^3}\frac{|\xi-v|^2}{2} f(t,x,\xi)d\xi.
\end{split}
\end{equation*}
For the traveling Maxwellian (\ref{traveling Maxwellian}), under transformation (\ref{A-0}), we have
\begin{equation}\label{M}
M=M_{[\br,\bv,\bt]}(\xi)=\frac{\br}{(2\pi r\bt)^{3/2}}\exp{(-\frac{|\xi-\bv|^2}{2r\bt})},
\end{equation}
where $r=2/3$, $\br(t,x)=\frac{2^{-3/2}}{R(t)^{3/2}}
e^{-\frac{h^2|x|^2}{2R(t)^2}}$, $\bv(t,x)=\frac{R'(t)}{R(t)}x$,
and $\bt(t,x)=\frac{3/4}{R(t)^2}$.

{\bf Theorem 1.2.} {\it For $\epsilon_0>0$ small, suppose that the
pulling speed $h\in(0,\epsilon_0^{1/2})$, and the initial data around the equilibrium $\mu$ satisfies
\begin{equation*}
    f_0=\mu+\mu^{1/2}\tilde{u}_0
\end{equation*}
with
\begin{equation}\label{decay of u}
\|\pb\tilde{u}_0\|_{\infty}<\epsilon_0.
\end{equation}
In addition, $f_0$ also satisfies the conservation laws
\bes
\int_{\Omega\times\R^3}(f_0-M)dyd\eta=0,
\ees
\bes
\int_{\Omega\times\R^3}(|\eta|^2+h^2|y|^2)(f_0-M)dyd\eta=0,
\ees
\bes
\int_{\Omega\times\R^3}(y\times\eta)(f_0-M)dyd\eta=0.
\ees
Then there exists a unique mild solution $f=M+M^{1/2}u\geq0$ to problem (\ref{eq: f}) together with (\ref{data: f}).
Moreover, scaling back to the original coordinates $(t,x,\xi)$ by (\ref{A-0}), we have the following decay estimates
\begin{equation}\label{decay of rho}
\frac{c_0}{R(t)^3}\leq\rho(t,x)\leq\frac{C_0}{R(t)^3},
\end{equation}
where $c_0$ and $C_0$ are two positive constants independent of $t$.

In addition, if the initial data $f_0(x,\xi)$ is continuous except on $\ga_{01}$, then $f(t,x,\xi)$ is continuous in $\{S\times\R^3\}\sm\ga_{\mathrm{sing}}$, where $\ga_{\mathrm{sing}}
=\{(t,x,\xi)\in\p S\times\R^3:\,t>0,\,|x|=R(t),\,\xi\cdot n_x=R'(t),\,|\xi|\geq h\}$,  $n_x$
stands for the unit outer normal at $x\in\p\O_t$, and
$\ga_{01}$ has been defined in Theorem 1.1.
}

\vskip 0.2 true cm

{\bf Remark 1.2.} {\it In \cite{Euler} and \cite{NS}, the movements of gases in an expanding ball $\Omega_t$
are globally described by the Euler
equations and the Navier-Stokes equations, respectively, and the authors have shown that
the vacuum will not appear in $\Omega_t$.}

\vskip 0.2 true cm

So far there exists extensive literature on the study of the Boltzmann equation. For the
Cauchy problem of the Boltzmann equation, under some special assumptions on the
initial data,
T.Carleman \cite{MR1555355} proved the local existence
of solutions for
the  spatially homogeneous case, while H.Grad \cite{MR0184507} and S.Ukai \cite{MR0363332} established the
local and global existence of solutions for the spatially inhomogeneous case. After that, by our knowledge,
the Boltzmann equation is mainly discussed in three different frameworks.
In $L^{\infty}$ framework, by the spectral analysis of the linearized Boltzmann equation, various initial value
problems and initial-boundary value
problems were considered (see \cite{MR839310}-\cite{MR714977} and so on). In $L^1$ framework, R.J.DiPerna
and P.L.Lions \cite{MR1014927} constructed the renormalized solutions, which were based on the velocity-averaging lemma and entropy dissipation
(see also \cite{MR1765272} and \cite{MR2116276}). In $L^2$ framework, based on macro-micro decomposition,
the authors in \cite{MR2043729} and \cite{MR2239361} developed energy methods for the Boltzmann equation,
and subsequently different wave patterns were considered by utilizing the energy methods
(e.g. \cite{MR2044894}, \cite{MR2221210}, \cite{MR2450610} and so on). In addition, the $L^2$ energy methods were generalized by Y.Guo and applied to various systems in \cite{MR1908664}-\cite{MR2000470} and \cite{MR2904573} respectively. On the other hand,
for the Boltzmann equation without angular cut-off, there also exist many results
(see \cite{MR2863853}, \cite{MR2793203}, \cite{MR2784329} and references therein).

For the initial-boundary problems of Boltzmann equation, Y.Guo developed the $L^2-L^{\infty}$ procedure in \cite{MR2679358} and solved the problems for all the four kinds of boundary conditions
(namely, the inflow, reverse-reflection, specular-reflection
and diffuse-reflection boundary conditions).
For the case of specular-reflection boundary condition, the domain was required to be strictly convex and analytic
in \cite{MR2679358} (when the domain is only strictly convex,
the analyticity assumption of boundary has been removed recently in \cite{2016arXiv160404342K}).
With respect to the regularities of solutions, the authors in \cite{MR3592757}
proved that the solution is $C^1$ away from the grazing set for the  Boltzmann equation
in a convex domain. While if the domain is non-convex, singularities of solution may propagate from the grazing set to the interior
of domain, see \cite{MR2855537} and \cite{MR3466841} for more details. The corresponding
results have been generalized to the Boltzmann equation with soft potential
and angular cut-off in \cite{MR3590377}.

To derive the $L^2$ decay in the $L^2-L^{\infty}$ procedure, the author  in \cite{MR2679358} applied the compactness method to a finite time interval. Since the
Boltzmann operator forms a semigroup, the long time $L^2$ decay of solutions was derived by iteration
in \cite{MR2679358}. Subsequently, by writing the linearized Boltzmann equation in a weak formulation,
a constructive method of $L^2$ estimate was constructed
in \cite{MR3085665} for the diffuse-reflection boundary condition. By choosing the test functions suitably, the $L^2$ estimates of the macro-components could be controlled by the micro-components, while the remaining boundary terms were
controlled well by the dissipation property of diffuse-reflection boundary.
The methods in \cite{MR2679358} and \cite{MR3085665} were generalized to the boundary condition
which was a linear combination of the
diffuse-reflection condition and specular-reflection condition in \cite{2015arXiv151101305B}.
However, the case of pure specular-reflection boundary condition was not considered
in these papers.

In this paper, we consider the Boltzmann equation in an expanding ball with pure specular-reflection boundary condition,
which leads to a new form (\ref{eq: f}) of Boltzmann equation whose coefficient depends on the time variable $\tau$.
In this case, the resulting linearized operator no longer forms a semigroup and we can not use the implicit
method directly to get the $L^2$ decay of solutions as in \cite{MR2679358}. Motivated by \cite{MR3085665},
we will apply the constructive method to the pure specular-reflection
boundary problem (\ref{eq: f}) although there is no dissipation property on the boundary.
Through choosing the Burnette functions as orthogonal bases in micro-components, we can reformulate the boundary integral
in a more delicate way and look for suitable test functions to handle the resulting boundary terms. As a byproduct, we give a constructive method to prove the $L^2$ decay for the Boltzmann equation in a bounded domain with specular-reflection boundary conditions. On the other hand, reverse-reflection boundary condition is ill-posed here because of the potential term in (\ref{eq: f}) (see Remark 2.1 below)

Here we point out that we have used the so-called ``traveling Maxwellian" in
(\ref{traveling Maxwellian}) to treat the Boltzmann equation (\ref{Boltzmann eq})
with (\ref{initial-boundary condition}) in the expanding ball $\O_t$.
The global Cauchy problem of the traveling Maxwellian was studied by R.Illner and M.Shinbrot in \cite{MR760333}.
For the extremely rarefied gases, the authors in \cite{MR760333} applied the iteration scheme in \cite{MR0475532}
to obtain a global mild solution of Boltzmann equation (one can also see Chapter 5 of \cite{MR1307620}).
In the present paper, since the gases are not extremely rarefied in the ball at the beginning
and lie in a bounded domain at any finite time,
we are required to give some different treatments from those in  \cite{MR760333}.

The paper is organized as follows: In Section 2, we discuss the properties of
backward trajectories of operator $\Lat$ and reformulate the Velocity Lemma
(see Lemma 2.4 below). In Section 3, we list some basic properties of Boltzmann operators which will be applied
later on. In Section 4, we establish the $L^2$-estimates of solutions to linear Boltzmann equations. In Section 5,
an explicit formula of solution to the transport equation is given, and subsequently
the $L^{\infty}$-estimate of solutions to a class of linear weighted Boltzmann equations
is established. In Section 6, at first, by the Duhamel's principle, one can write
out the implicit expression of the solution to full Boltzmann equation. Based on this, by iteration,
we derive the existence and uniqueness of the solution to problem (\ref{eq: f}) with (\ref{data: f}).
And then the proofs of Theorem 1.1 and Theorem 1.2 are completed.

\textbf{Notations}: In the following sections, to simplify the notations, we denote by $\|u\|=\|u\|_{L^2(\O\times\R^3)}$, $\|u\|_{\infty}=\|u\|_{L^{\infty}(\O\times\R^3)}$ and $\|u\|_{\nu}=\|\sqrt{\nu}u\|_{L^2(\O\times\R^3)}$,
where $\O=\{y=(y_1, y_2,y_3)\in\Bbb R^3: |y|<1\}$ and $\nu=
\nu(\eta)=\sqrt{2\pi}[(|\eta|+|\eta|^{-1})
\int_0^{|\eta|}e^{-s^2}ds+e^{-|\eta|^2}]$. We say $A\lesssim B$ if $A\leq CB$ holds for some positive constant $C$
independent of the quantities $A$ and $B$.

\section{Backward trajectories of operator $\Lat$}
Recall that $\Lat=\p_{\tau}+\eta\cdot\nabla_{y}-h^2y\cdot\nabla_{\eta}$ in \eqref{A-A}.
We now consider  the trajectory $(Y,H)(\tau)=(Y,H)(\tau;\tau_0,y_0,\eta_0)$  of $\Lat$ in the whole phase space $\R^3_{y}\times\R^3_{\eta}$
starting from the point $(y_0,\eta_0)$ and time $\tau_0$:
\begin{equation*}
\left\{ \enspace
\begin{aligned}
&\frac{dY}{d\tau}=H,\\
&\frac{dH}{d\tau}=-h^2Y,\\
&(Y,H)|_{\tau=\tau_0}=(y_0,\eta_0).
\end{aligned}
\right.
\end{equation*}
Direct computation yields
\begin{equation}
\left\{ \enspace
\begin{aligned}
Y(\tau)&=y_0\cos[h(\tau-\tau_0)]
+\frac{1}{h}\eta_0\sin[h(\tau-\tau_0)],\\
H(\tau)&=-hy_0\sin[h(\tau-\tau_0)]
+\eta_0\cos[h(\tau-\tau_0)].
\end{aligned}
\right.\label{eq: Y and H}
\end{equation}
This means that $Y(\tau)$ is an ellipse in $\R^3_y$ with the origin $O$ as its center. Since the related
time interval is $(0,2\pi/h)$, the trajectory only takes a quarter of the ellipse.
Moreover, we have the conservations of energy and angular momentum for $(H, Y)$:
\bes
    \frac{d}{d\tau}(|H|^2+h^2|Y|^2)=0
\ees
and
\bes
\frac{d}{d\tau}(Y\times H)=0.
\ees
From this, we can define the numbers $e_0$ and $m_0$ as follows
\be\label{eq: e0}
e_0:=|\eta_0|^2+h^2|y_0|^2=|H(\tau)|^2+h^2|Y(\tau)|^2,
\ee
\be\label{eq: m0}
m_0:=|\eta_0\times y_0|^2=|Y(\tau)\times H(\tau)|^2.
\ee
Denote by $l_{\max}$ and $l_{\min}$ the major and minor semi-axis of ellipse
$Y(\tau)$, respectively. It is easy to know that
\bes
|l_{\max}|^2=\left(e_0+\sqrt{e_0^2-4h^2m_0}\right)/2h^2,
\ees
\bes
|l_{\min}|^2=\left(e_0-\sqrt{e_0^2-4h^2m_0}\right)/2h^2.
\ees

\vskip 0.2 true cm

Next we study the properties of operator $\Lat$ in the domain $\bar{\O}\times\R^3_{\eta}$,
which will be divided into two cases of $|y_0|<1$ and $|y_0|=1$.

\subsection{Backward trajectory $Y(\tau)$ in the interior of $\O$}
For $(y_0,\eta_0)\in\O\times\R^3$, the backward trajectory $Y(\tau)$ may hit the boundary $\p\O$ and change its velocity,
then travels along another ellipse. This can be precisely stated as follows

\vskip 0.2 true cm

{\bf Lemma 2.1.} {\it
For $(y_0,\eta_0)\in\O\times\R^3$, we have

(a) If  $l_{\min}<1<l_{\max}$, i.e., $e_0-m_0>h^2$, then the backward trajectory  $Y(\tau)$
will hit the boundary $\p\O$ at some point, then reflect specularly and travel along another ellipse.
In particular, this is the case when $|\eta_0|\geq2h$;

(b) If $l_{\min}<l_{\max}=1$, i.e., $e_0-m_0=h^2$, then the backward trajectory $Y(\tau)$
will graze the boundary $\p\O$ at some point and travel along the same ellipse;

(c) If $l_{\min}\leq l_{\max}<1$, i.e., $e_0-m_0<h^2$, then the backward trajectory $Y(\tau)$
remains in the interior of $\O$ and travels along the same ellipse.
}

{\bf Proof.} Note that $l_{\min}<1$ holds for $y_0\in\O$. Thus it only suffices
to consider the relation between $l_{\max}$ and $1$. In fact, Lemma 2.1 holds
by direct verifications and observation.
\qquad \qquad \qquad \qquad \qquad  $\square$

\vskip 0.2 true cm

{\bf Definition 2.2.} {\it  For $(y_0,\eta_0)\in\O\times\R^3$ satisfying $e_0-m_0>h^2$
(corresponding to Case (a) of Lemma 2.1), we define
  \be\label{def: backward exet}
  \begin{aligned}
    \tau_b(y_0,\eta_0):= & \sup\{\tau>0:\,Y(\tau_0-s;\tau_0,y_0,\eta_0)\in\O \text{ for all } 0<s<\tau\},\\
    y_b(y_0,\eta_0):= & Y(\tau_0-\tau_b;\tau_0,y_0,\eta_0), \\
    \eta_b(y_0,\eta_0):= & H(\tau_0-\tau_b;\tau_0,y_0,\eta_0).
  \end{aligned}
  \ee
}
Here we point out that $(\tau_0-\tau_b,y_b,\eta_b)$ is just the point when the backward trajectory
$H(\tau;\tau_0,y_0,\eta_0)$ starting from $(\tau_0,y_0,\eta_0)$ first hits the boundary $\p\O$.
In order to define the backward trajectory piece by piece, we introduce the following notation.

\vskip 0.2 true cm

{\bf Definition 2.3.} {\it Let $(y_0,\eta_0)\in\O\times\R^3$ satisfy $e_0-m_0>h^2$, for $k\geq 0$, we define
  \begin{align*}
    \tau_{k+1}:= & \tau_k-\tau_b(y_k,\eta_k), \\
    y_{k+1}:= & y_b(y_k,\eta_k), \\
    \eta_{k+1}:= & \eta_b(y_k,\eta_k)-2[y_{k+1}\cdot\eta_b(y_k,\eta_k)]y_{k+1}.
  \end{align*}}
In this case,  the backward trajectory $(y(\tau),\eta(\tau))$ starting from $(\tau_0,y_0,\eta_0)$ can be expressed as
\be\label{eq: y eta backward}
 \ba
  y(\tau):= & \sum_{k=0}^{\infty}1_{(\tau_{k+1},\tau_k]}(\tau)
    Y(\tau;\tau_k,y_k,\eta_k),\\
  \eta(\tau):= & \sum_{k=0}^{\infty}1_{(\tau_{k+1},\tau_k]}(\tau)
  H(\tau;\tau_k,y_k,\eta_k),
 \ea
\ee
where $1_{(\tau_{k+1},\tau_k]}(\tau)$ stands for the characteristic function of interval $(\tau_{k+1},\tau_k]$.

For Case (b) of Lemma 2.1, the representations in \eqref{eq: y eta backward} for the backward trajectory are still plausible.
The only difference from Case (a) is that the trajectory grazes the boundary $\p\O$. This means that there is actually no change of velocity
$u$ at the grazing point due to the specular-reflection boundary condition.

For Case (c) of Lemma 2.1, the backward exit time $\tau_b(y_0,\eta_0)$ cannot be defined since the trajectory remains
in the interior of domain $\O$. But we  still use the representation \eqref{eq: y eta backward}
with only $k=0$.

To derive the $L^{\infty}$ decay of solutions to linearized Boltzmann equations,
for small $\ka>0$, large $N>0$,
we need to define the following set
\be\label{def: set A}
  \AehN:=\{(y_0,\eta_0)\in\O\times\R^3:\,
    e_0-m_0\geq h^2+\ka^2,\,2h\leq|\eta_0|\leq2N\}.
\ee

{\bf Lemma 2.4 (Velocity Lemma)}  {\it
For $(\tau_0,y_0,\eta_0)\in(0,\pi/2h)\times\AehN$, denote the backward trajectory
$(y(\tau),\eta(\tau))$ by (\ref{eq: y eta backward}), then
we have

(a) The time interval $\Delta\tau$ between two adjacent reflections point is
\be\label{eq: Delta tau}
\Delta\tau=\frac{1}{h}
\arccos\left(\frac{e_0-2h^2}{\sqrt{e_0^2-4h^2m_0}}\right),
\ee
and admits the upper and lower bounds as
\be\label{estimates: Delta tau}
\frac{\ka}{2N^2}\leq\Delta\tau\leq\frac{\pi}{2h}.
\ee

(b)
For $k>\frac{\pi N^2}{h\ka}$, we have $\tau_k<0$, which means
that the summation of $k$ in (\ref{eq: y eta backward}) is finite.

(c)
For $k\leq\frac{\pi N^2}{h\ka}$, set $\de=\frac{h^2\ka}{N^2}$. When $\tau\in(\tau_k+\de,\tau_{k-1}-\de)$, we have
that for all $\eta'\in\R^3$ satisfying $2h\leq |\eta'|\leq 2N$,
\be\label{estimates: Y1 H1}
(y(\tau),\eta')\in\AehhN.
\ee

(d)
The measure of the set $(0,\tau_0)\backslash\cup_k(\tau_k-\de,\tau_k+\de)$ is less than $2\pi h$.

(e)
$\tau_k,\,y_k$ and $\eta_k$ are analytic functions of $(\tau_0,y_0,\eta_0)$.
}

{\bf Proof.} (a) For $(\tau_0,y_0,\eta_0)\in(0,\pi/2h)\times\AehN$, from (\ref{eq: Y and H}) we have
\bes
Y(\tau)=y_0\cos[h(\tau-\tau_0)]+\frac{\eta_0}{h} \sin[h(\tau-\tau_0)].
\ees
Then we get
\be\label{eq: Y(tau)2}
|Y(\tau)|^2=\frac{1}{2h^2}
\left(e_0-\sqrt{e_0^2-4h^2m_0}\cos[2h(\tau-\tau_0)+\theta]\right),
\ee
where
\bes
\cos\theta=\frac{|\eta_0|^2-h^2|y_0|^2}{\sqrt{e_0^2-4h^2m_0}}>0, \,\, \sin\theta=\frac{2hy_0\eta_0}{\sqrt{e_0^2-4h^2m_0}}.
\ees
Let $|Y(\tau)|=1$ in (\ref{eq: Y(tau)2}) yield
\bes
\cos[2h(\tau-\tau_0)+\theta]=\frac{e_0-2h^2}{\sqrt{e_0^2-4h^2m_0}}>0.
\ees
This derives
\bes
\tau_k=\tau_0\pm\frac{1}{2h}
\arccos\left(\frac{e_0-2h^2}{\sqrt{e_0^2-4h^2m_0}}\right)
-\frac{2k\pi+\theta}{2h}.
\ees
Hence,  the backward exit time $\tau_b$ is
\be\label{eq: tau b}
\tau_b(y_0,\eta_0)=\frac{1}{2h}
\arccos\left(\frac{e_0-2h^2}{\sqrt{e_0^2-4h^2m_0}}\right)
-\frac{1}{2h}
\arccos\left(\frac{e_0-2h^2|y_0|^2}{\sqrt{e_0^2-4h^2m_0}}\right),
\ee
and the time interval between two adjacent reflection is
\bes
\Delta\tau=\frac{1}{h}
\arccos\left(\frac{e_0-2h^2}{\sqrt{e_0^2-4h^2m_0}}\right).
\ees
For the upper bound of $\Delta\tau$, by $e_0^2\geq\eta_0^2\geq 4 h^2$, we have
\bes
\Delta\tau\leq\frac{\pi}{2h}.
\ees
In addition, by $e_0-m_0\geq h^2+\ka^2$, we have
\bes
\frac{e_0-2h^2}{\sqrt{e_0^2-4h^2m_0}}\leq
\frac{e_0-2h^2}{\sqrt{(e_0-2h^2)^2+4h^2\ka^2}}.
\ees
This yields
\bes\ba
\Delta\tau\geq & \frac{1}{h}
\arcsin\left(\frac{2h\ka}{\sqrt{(e_0-2h^2)^2+4h^2\ka^2}}\right) \\
\geq & \frac{\ka}{2N^2}.
\ea\ees

(b) Since $\tau_0\in(0,\frac{\pi}{2h})$ and $\Delta\tau\geq\frac{\ka}{2N^2}$ from (a), the number
of reflections along the backward trajectory is less than $\frac{\pi N^2}{h\ka}$. Hence we have
\bes
\tau_k<0\quad \text{for $k\geq\frac{\pi N^2}{h\ka}$}.
\ees

(c) Since $\delta=\frac{h^2\ka}{N^2}<\Delta\tau$, we have
\bes
(2h(\tau_k-\tau_0)+\theta,2h(\tau_{k-1}-\tau_0)-\theta)
\in(2k\pi-\frac{\pi}{2},2k\pi+\frac{\pi}{2}).
\ees
From (\ref{eq: Y(tau)2}), $|y(\tau)|^2$ is a convex function of $\tau\in(\tau_k+\de,\tau_{k-1}-\de)$. Thus,
\be\label{estimate: y by l}
\begin{aligned}
|y(\tau)|^2\leq & |y(\tau_k+\de)|^2 \\
\leq & (1-\frac{2\de}{\Delta\tau})|y(\tau_k)|^2
+\frac{2\de}{\Delta\tau}|y(\tau_k+\frac{\Delta\tau}{2})|^2\\
\leq & (1-\frac{2\de}{\Delta\tau})+\frac{2\de}{\Delta\tau}|l_{\min}|^2\\
\leq & 1-\frac{2\de}{\Delta\tau}(1-|l_{\min}|^2).
\end{aligned}
\ee
It follows from $e_0-m_0\geq h^2+\ka^2$ that
\bes
\begin{aligned}
\frac{1-|l_{\min}|^2}{\Delta\tau} &
=\frac{\sqrt{e_0^2-4h^2m_0}}{2h}\cdot
\frac{1-\frac{e_0-2h^2}{\sqrt{e_0^2-4h^2m_0}}}
{\arccos\left(\frac{e_0-2h^2}{\sqrt{e_0^2-4h^2m_0}}\right)} \\
& \geq\frac{\sqrt{e_0^2-4h^2m_0}}{2h}\cdot
\frac{1-\frac{e_0-2h^2}{\sqrt{(e_0-2h^2)^2+4h^2\ka^2}}}
{\arccos\left(\frac{e_0-2h^2}{\sqrt{(e_0-2h^2)^2+4h^2\ka^2}}\right)}\\
& \geq\frac{1}{2h}\cdot
\frac{\sqrt{(e_0-2h^2)^2+4h^2\ka^2}-(e_0-2h^2)}
{\arcsin\left(\frac{2h\ka}{\sqrt{(e_0-2h^2)^2+4h^2\ka^2}}\right)}\\
& \geq \frac{\ka}{2}.
\end{aligned}
\ees
This, together with (\ref{estimate: y by l}), yields
\bes
|y(\tau)|^2\leq 1-\de\ka=1-\frac{h^2\ka^2}{N^2}.
\ees
Note that for any $\eta'\in\R^3$ satisfying $2h\leq|\eta'|\leq 2N$, we have
\bes
\begin{aligned}
e'-m'= & \{|\eta'|^2+h^2|y(\tau)|^2\}
-|\eta'\times y(\tau)|^2 \\
\geq & h^2+(|\eta'|^2-h^2)(1-|y(\tau)|^2)\\
\geq & h^2+\frac{h^4\ka^2}{N^2}.
\end{aligned}
\ees
Then for $\tau\in(\tau_k+\de,\tau_{k-1}-\de)$, $2h\leq|\eta'|\leq 2N$, we conclude that
\bes
(y(\tau),\eta')\in\AehhN.
\ees

(d) The measure of the set $(0,\tau_0)/\cup_k(\tau_{k+1}+\de,\tau_k-\de)=\cup_k[\tau_k-\de,\tau_k+\de]$ is less than
\bes
2\de\cdot \frac{\pi N^2}{h\ka}\leq 2\pi h.
\ees

(e) To prove the analyticity of $\tau_k,\,y_k$ and $\eta_k$ with respect
to the variable $(\tau_0,y_0,\eta_0)$, for $(y_0,\eta_0)\in\AehN$, we only need to
study the functions $\tau_b$, $y_b$
and $\eta_b$ defined in (\ref{def: backward exet}). Since $|y_b|=1$, $y_b\cdot\eta_b<0$ and
\bes\ba
\ka^2\leq & e_0-m_0-h^2 \\
= & (|\eta_b|^2-h^2)(1-|y_b|^2)+(y_b\cdot\eta_b)^2\\
= & (y_b\cdot\eta_b)^2,
\ea\ees
we have
\bes
y_b\cdot\eta_b\leq -\ka.
\ees
When solving $|y(-\tau_b)|^2=1$, by the fact that
\bes
\frac{\p}{\p s}|y(-s)|^2\bigg\arrowvert_{s=\tau_b}
=-y_b\cdot\eta_b\geq \ka>0,
\ees
then we can see that $\tau_b$ is locally solvable. Similarly, direct computation yields
that the derivatives
of $|y(-\tau_b)|^2$ with respect to variables $y$ and $\eta$ are
\bes
(\nabla_{y_0}\tau_b,\nabla_{\eta_0}\tau_b)=
\left(-\frac{\cos(h\tau_b)}{(y_b\cdot\eta_b)}y_b,
 -\frac{\sin(h\tau_b)}{h(y_b\cdot\eta_b)}y_b\right).
\ees
Note that $|y(\tau)|^2$ is analytic. Then by implicit function theorem,
we know that   $\tau_b$, $y_b$ and $\eta_b$ are analytic with respect
to the variable $(\tau_0,y_0,\eta_0)$.\qquad \qquad \qquad \qquad \qquad \qquad
\qquad \qquad \qquad \qquad \qquad \qquad \quad $\square$

\subsection{Backward trajectory near boundary $\p\O$}
For $(y_0,\eta_0)\in\p\O\times\R^3$, the property of backward trajectory is more subtle. We denote
the phrase boundary as $\ga=\p\O\times\R^3$, and split $\ga$ into the outgoing boundary
$\ga_+$, the incoming boundary $\ga_-$, and the grazing set $\ga_0$ as follows:
\begin{align*}
\ga_+&=\{(y,\eta)\in\p S_0\times\R^3:\eta\cdot n_y>0\},\\
\ga_-&=\{(y,\eta)\in\p S_0\times\R^3:\eta\cdot n_y<0\},\\
\ga_0&=\{(y,\eta)\in\p S_0\times\R^3:\eta\cdot n_y=0\}.
\end{align*}
Compared with \cite{MR2679358} where $\ga_0$ is a singular set, in the present paper, only some part of $\ga_0$ is singular.
Since the potential force $-h^2y$ in $\Lat$ is pointing to the center of the unit ball, particles on part of $\ga_0$
will depart from the boundary and go to the interior of the ball. We should further split $\ga_0$ into non-singular
set $\ga_{00}$, and singular set $\ga_{01}$ as follows:
\be\ba\label{eq: singular set gamma}
\ga_{00}= & \{(y,\eta)\in\ga_0: |\eta|<h\},\\
\ga_{01}= & \{(y,\eta)\in\ga_0: |\eta|\geq h\}.
\ea\ee
The action of non-singular set $\ga_{00}$ is similar to that of the interior domain
of $\O\times\R^3$ while the singular set $\ga_{01}$ acts as the singular grazing set. More precisely, we have the following conclusion

\vskip 0.2 true cm

{\bf Lemma 2.5.} {\it The backward trajectory $(y,\eta)(\tau;\tau_0,y_0,\eta_0)$ is continuous for all
\bes
(\tau_0,y_0,\eta_0)\in(0,\pi/2h)\times
\{(\bar{\O}\times\R^3)\backslash \ga_{01}\}.
\ees
}

{\bf Proof.} To prove Lemma 2.5, it only suffices  to study the situation around the grazing set $\ga_{00}$ where there is at most one
collision with the boundary. In this case, we require to consider three classes of points, $(y_j,\eta_j)$, $j=1,2,3$, representing
the three cases discussed in Lemma 2.1 (see Figure \ref{fig:2} below) respectively.

\begin{figure}[htbp]
  \centering\includegraphics[width=6.5cm,height=6.5cm]{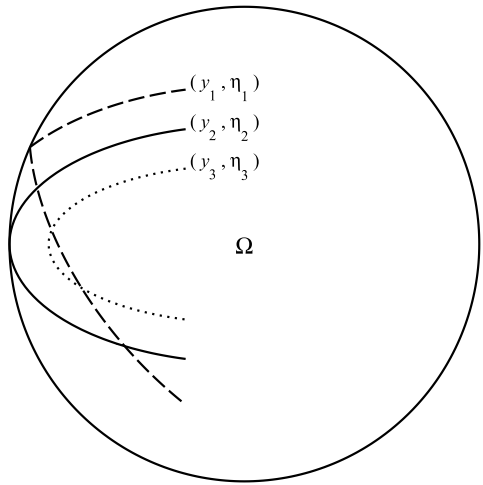}
  \caption{\bf Continuity for specular reflection boundary}\label{fig:2}
\end{figure}
For $j=1,2,3$, set
\bes
e_j=|\eta_j|^2+h^2|y_j|^2,
\ees
\bes
m_j=|\eta_j\times y_j|^2.
\ees
Then
\bes
e_1-m_1>h^2,\,\,e_2-m_2=h^2,\,\,e_3-m_3<h^2.
\ees

In the first two cases of Lemma 2.1, when $|(\tau_1,y_1,\eta_1)-(\tau_2,y_2,\eta_2)|<\ve$ sufficiently small,
it follows from \eqref{eq: tau b}  that for $j=1,2$,
\bes
\tau_{b_j}=:\tau_b(y_j,\eta_j)=\frac{1}{2h}
\arccos\left(\frac{e_j-2h^2}{\sqrt{e_j^2-4h^2m_j}}\right)
-\frac{1}{2h}
\arccos\left(\frac{e_j-2h^2|y_j|^2}{\sqrt{e_j^2-4h^2m_j}}\right).
\ees
Together with the fact
\bes\ba
e_j^2-4h^2m_j= & 2h^2(e_j-m_j-h^2)+(|\eta_j|^2-h^2)^2 \\
& +2h^2|y_j\cdot\eta_j|^2+h^4(1-|y_j|^2)^2 \\
\geq & h^4(1-|y_j|^2)^2,
\ea\ees
this yields that $\nabla_{y,\eta}\tau_b$ is locally finite around $(y_j,\eta_j)$
and $|\tau_{b1}-\tau_{b2}|\leq C\ve$ holds. After specular reflection, it is easy
to know $|y_1(\tau)-y_2(\tau)|\leq C\ve$.

In the last two case, when $|(\tau_2,y_2,\eta_2)-(\tau_3,y_3,\eta_3)|<\ve$ sufficiently small,
we have $y_j(\tau)=Y(\tau;\tau_j,y_j,\eta_j)$, $j=2,3$. This derives $|y_2(\tau)-y_3(\tau)|\leq C\ve$.
\qquad \qquad \qquad \qquad \qquad \qquad \qquad \quad $\square$

\begin{figure}[htbp]
\centering\includegraphics[width=6.5cm,height=6.5cm]{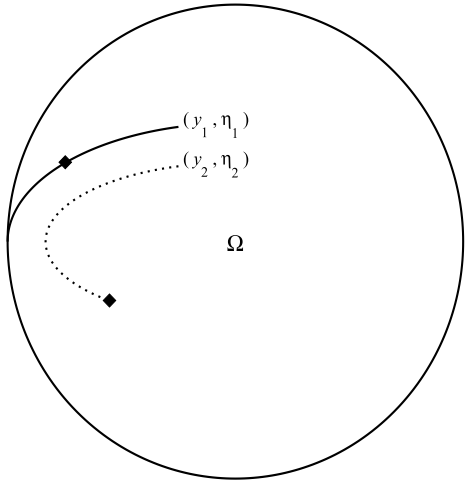}
\caption{\bf Discontinuity for reverse reflection boundary}\label{fig:3}
\end{figure}

{\bf Remark 2.1.} {\it When considering the same problem for the reverse-reflection boundary, we cannot
get any continuity result as in Lemma 2.5. As Figure \ref{fig:3} shows,
if we trace back along the trajectories of $(y_1,\eta_1)$ and $(y_2,\eta_2)$,
then the backward trajectory of $(y_2,\eta_2)$
stays on the same ellipse. But the backward trajectory of $(y_1,\eta_1)$ hits the boundary $\p\O$ and
then reflects reversely, which leads to the fact that these two trajectories of $(y_1,\eta_1)$ and $(y_2,\eta_2)$
can no longer stay close to each other.}

\section{Basic properties of the Boltzmann operators}
Before establishing the $L^2$ estimate of solutions to the linear Boltzmann equation, we list some properties of
the Boltzmann operators. The proofs are elementary and can be found in \cite{MR0156656} and Chapter 7 of \cite{MR1307620}.

The linearized collision operator $L$, given by (\ref{L}), is a self-adjoint nonnegative operator in $L^2$.
The null space $\mathcal{N}$ of $L$ is spanned by
\bes
\mu^{1/2},\,\eta^j\mu^{1/2} (j=1,2,3),\,\frac{|\eta|^2}{2}\mu^{1/2},
\ees
which are known as the collision invariants. After normalization, we write
\bes
\begin{aligned}
\chi_0= & \mu^{1/2}, \\
\chi_j= & \eta^j\mu^{1/2},\, j=1,2,3, \\
\chi_4= & \frac{|\eta|^2-3}{\sqrt{6}}\mu^{1/2}.
\end{aligned}
\ees
In addition, we denote the projection to the null space $\mathcal{N}$ by $P$ as follows
\be\label{def: Pu}
Pu:=\sum_{k=0}^{4}<u,\chi_k>\chi_k,
\ee
where $<g,h>=\int_{\R^3}ghd\eta$. As in \cite{MR2239361},
we shall use the following Burnette functions of the space $\mathcal{N}^{\perp}$:
\bes
\begin{aligned}
A_j(\eta)= & \eta^j\frac{|\eta|^2-5}{\sqrt{10}}\mu^{1/2}, &
& \text{for $j=1,2,3,$} \\
B_{kl}(\eta)= & \left(\eta^k\eta^l-\frac{\de_{kl}}{3}|\eta|^2\right)\mu^{1/2},&
& \text{for $k,l=1,2,3$},
\end{aligned}
\ees
where $\de_{kl}=1$ if $k=l$ and $\de_{kl}=0$ if $k\neq l$.
Direct verification yields

\vskip 0.2 true cm
{\bf Lemma 3.1.} {\it For $i,j,k,l=1,2,3$,

(a) $PA_j=0$ and $PB_{kl}=0$.

(b) $<A_j,A_i>=\de_{ji}$ and $<A_j,B_{kl}>=0$.

(c) $<B_{ij},B_{kl}>=\de_{ik}\de_{jl}+\de_{il}\de_{jk}
-\frac{2}{3}\de_{ij}\de_{kl}$.
}

\vskip 0.2 true cm

On the other hand, the operator $L$ can be split as $-\nu(\eta)I+K$, where
\begin{equation}\label{eq: nu}
\nu(\eta)=\nu(|\eta|)=\sqrt{2\pi}\left[(|\eta|+|\eta|^{-1})
\int_0^{|\eta|}e^{-s^2}ds+e^{-|\eta|^2}\right],
\end{equation}
and $K$ is an integral operator with the symmetric kernel function $k(\eta,\eta_*)$ given by
\be\ba\label{eq: k}
k(\eta,\eta_*)= & \sqrt{2\pi}|\eta_*-\eta|^{-1}
\exp\left(-\frac{(|\eta_*|^2-|\eta|^2)^2}{8|\eta_*-\eta|^2}
-\frac{|\eta_*-\eta|^2}{8}\right) \\
& -\frac{1}{2}|\eta_*-\eta|
\exp\left(-\frac{|\eta_*|^2+|\eta|^2}{4}\right).
\ea\ee
Obviously, $\nu(\eta)$ satisfies
\begin{equation}\label{bound for nu}
0<\nu_0\leq\nu(\eta)\leq\nu_1(1+|\eta|),
\end{equation}
where $\nu_0$ and $\nu_1$ are certain positive constants. And the kernel function $k(\eta,\eta_*)$ admits the following
properties:

\vskip 0.2 true cm

{\bf Lemma 3.2.} {\it $k(\eta,\eta_*)$ is integrable and square integrable
with respect to the variable $\eta_*\in\Bbb R^3$, moreover, the integral is bounded by a
positive constant independent of $\eta$.
}

\vskip 0.2 true cm

{\bf Proof.}
Following the expression of $k(\eta,\eta_*)$ in (\ref{eq: k}), we have
\begin{equation}\label{estimate: k}
|k(\eta,\eta_*)|\leq C \left(|\eta_*-\eta|^{-1}
+|\eta_*-\eta|\right)\exp{\left(-\frac{1}{8}|\eta_*-\eta|^2\right)},
\end{equation}
where $C>0$ is a constant independent of $\eta$.
Then
\begin{equation*}
\int_{\R^3}|k(\eta,\eta_*)|d\eta_*\leq C\int_{\R^3}\left(|\eta_*|^{-1}
+|\eta_*|\right)\exp{\left(-\frac{1}{8}|\eta_*|^2\right)}d\eta_*
\leq C.
\end{equation*}
Similar result holds for the square norm of $k(\eta,\eta_*)$ with respect to
variable $\eta_*$.
\qquad \qquad \qquad \qquad \quad $\square$

Let the weight function $\phi(y,\eta)=(1+|\eta|^2+h^2|y|^2)^{\beta/2}$ with $\beta>3/2$. Then we have

\vskip 0.2 true cm

{\bf Lemma 3.3.} {\it The operator $K_{\phi}$ is a bounded operator in $L^{\infty}(\Bbb R^3)$
space, where $K_{\phi}w=\phi K\left(\phi^{-1}w\right)$.}

\vskip 0.2 true cm

{\bf Proof.} The kernel of the operator $K_{\phi}$ is
\begin{equation*}
k_{\phi}(\eta,\eta_*)=k(\eta,\eta_*)\frac{\phi(y,\eta)}{\phi(y,\eta_*)}.
\end{equation*}
Note that for fixed $y\in\O$,
\begin{equation*}
\phi(y,\eta)\leq\phi(y,\eta_*)\phi(y,\eta_*-\eta).
\end{equation*}
Then we have
\begin{equation*}
|\kp(\eta,\eta_*)|\leq C \phi(y,\eta_*-\eta)\left(|\eta_*-\eta|^{-1}
+|\eta_*-\eta|\right)\exp{\left(-\frac{1}{8}|\eta_*-\eta|^2\right)}.
\end{equation*}
Hence $K_{\phi}$ is a bounded operator in $L^{\infty}$.\qquad \qquad \qquad \qquad \qquad \qquad
\qquad \qquad \qquad \qquad \qquad $\square$

\vskip 0.2 true cm

It is well-known that the operator $L$ satisfies
\begin{equation}\label{negetivity of L}
  <Lu,u>\leq -\sigma\|(I-P)u\|_{\nu}^2,
\end{equation}
where $\|h\|_{\nu}=\|\nu(\eta)^{1/2}h\|$. As for the bilinear operator $\Gamma(g,h)$, we have

\vskip 0.1 true cm

{\bf Lemma 3.4.} {\it
There exists a constant $C>0$ independent of $\eta$ such that
\begin{equation*}
|\Gamma(g,h)(\eta)|\leq C\nu(\eta)\|g\|_{\infty}\|h\|_{\infty}.
\end{equation*}
}

{\bf Proof.}
From the definitions in (\ref{Gamma}) and (\ref{Q(f,g)}), we can split $\Gamma(g,h)$ as
\begin{equation*}
\Gamma(g,h)=\Gamma_1+\Gamma_2-\Gamma_3-\Gamma_4.
\end{equation*}
Note $\mu\mu_*=\mu'\mu'_*$, we then have
\begin{align*}
|\Gamma_1|&\leq 2\pi\mu^{-1/2}\int_{\R^3}|\eta_*-\eta|(\mu^{1/2}g)'(\mu^{1/2}h)'_*d\eta_*\\
&\leq 2\pi\|g\|_{\infty}\|h\|_{\infty}\int_{\R^3}|\eta_*-\eta|\mu_*^{1/2}d\eta_*\\
&\leq C\nu(\eta)\|g\|_{\infty}\|h\|_{\infty},
\end{align*}
where in the last inequality we have used the property of $\nu(\eta)$ in (\ref{bound for nu}). In addition,
\begin{align*}
|\Gamma_3|&\leq 2\pi\mu^{-1/2}\int_{\R^3}|\eta_*-\eta|(\mu^{1/2}g)(\mu^{1/2}h)_*d\eta_*\\
&\leq 2\pi\|g\|_{\infty}\|h\|_{\infty}\int_{\R^3}|\eta_*-\eta|\mu_*^{1/2}d\eta_*\\
&\leq C\nu(\eta)\|g\|_{\infty}\|h\|_{\infty}.
\end{align*}
Similar estimates hold for $\Gamma_2$ and $\Gamma_4$. Thus, Lemma 3.4 is proved.\qquad
\qquad \qquad \qquad \qquad \qquad $\square$

\vskip 0.2 true cm

Lemma 3.4 together with the fact of $\phi(\eta)\leq\phi(\eta')\phi(\eta'_*)$, yields

\vskip 0.1 true cm

{\bf Lemma 3.5.} {\it There exists a constant $C>0$ independent of $\eta$ such that
\begin{equation*}
|\Gamma_{\phi}(g,h)(\eta)|\leq C\nu(\eta)\|g\|_{\infty}\|h\|_{\infty},
\end{equation*}
where $\Gamma_{\phi}(g,h)=\phi\Gamma\left(\phi^{-1}g,\phi^{-1}h\right)$.
}

\section{$L^2$-estimate of solutios to the linear Boltzmann equation}
Consider the following linear Boltzmann equation
\begin{equation}\label{eq: linear u with g}
\pt u+\eta\cdot\nabla_yu-h^2y\cdot\nabla_{\eta}u
=\mut\cos^2(h\tau)Lu+\cos^2(h\tau)g,
\end{equation}
where $g=g(\tau,y,\eta)$ is a smooth function. The initial-boundary condition of $u$ is given by
\begin{equation}\label{data for u in L2}
\begin{cases}
 u(0,y,\eta)=u_0(y,\eta),&
 \text{$y\in\O,\,\eta\in\R^3$,}\\
 u(\tau,y,\eta)=u(\tau,y,\eta-2(\eta\cdot n_y)n_y),&
 \text{$y\in\p\O$, $\eta\cdot n_y<0$.}
\end{cases}
\end{equation}
In addition, we assume that $g$ also satisfies the following condition to assure the conservation of mass, energy and angular momentum:
\begin{equation}\label{mass conservation for g}
Pg=0,
\end{equation}
where the operator $P$ is defined in \eqref{def: Pu}.
Set
\bes
\al(\tau):=\int_{0}^{\tau}\cos^2(hs)ds.
\ees
We will prove the following $L^2$ decay estimates of solution $u$ to \eqref{eq: linear u with g}.

\vskip 0.2 true cm

{\bf Proposition 4.1.}  {\it
Let $u\in L^2$ be the weak solution of (\ref{eq: linear u with g}) with initial-boundary value condition (\ref{data for u in L2}), then there exists a constant $\la>0$ such that for $\tau\in(0,\pi/2h)$,
\begin{equation}\label{estimates: u linear}
\|u(\tau)\|^2\lesssim e^{-2\la\al(\tau)}\{\|u(0)\|^2+\int_{0}^{\tau}
e^{2\la\al(s)}\|g(s)\|^2ds\}.
\end{equation}
}

The rest of this section is devoted to the proof of Proposition 4.1.
In terms of the macro-micro decomposition, set
\bes\ba
a(\tau,y)= & <u(\tau,y,\eta),\chi_0>\mut^{1/2}, \\
b^j(\tau,y)= & <u(\tau,y,\eta),\chi_j>\mut^{1/2}, \,\,j=1,2,3,\\
c(\tau,y)= & <u(\tau,y,\eta),\chi_0>\mut^{1/2}+q(\tau,y)\mut^{1/2}, \\
d(\tau,y,\eta)= & (I-P)u(\tau,y,\eta)\mut^{1/2},
\ea\ees
where $\mut=e^{-h^2|y|^2/2}$ and $q=\frac{h^2|y|^2}{\sqrt{6}}a$.
Let
\be\label{eq: u of abc}
u=\big\{a\chi_0+\sum_{j=1}^{3}b^j\chi_j+(c-q)\chi_4\big\}\tmn+d\tmn.
\ee
It follows from the conservation laws of mass and energy in (\ref{law of mass}) and (\ref{law of energy}) that
\be\label{eq: conserve a}
\int_{\O}ady=0
\ee
and
\be\label{eq: conserve c}
\int_{\O}cdy=0.
\ee
In addition,
\begin{equation}\label{estimats: q}
\|q\|^2\lesssim h^2\|a\|^2.
\end{equation}

\vskip 0.2 true cm

Next, we derive the estimates of $a, b=(b^1,b^2,b^3), c$ in terms of $d$ and $g$.
Rewrite the linear Boltzmann equation (\ref{eq: linear u with g}) in the following weak formulation:
\begin{equation}\label{eq: weak form for linear:1st}
\begin{split}
&\iOR\left\{\psi u(\tau)-\psi u(s)\right\}dyd\eta-\int_{s}^{\tau}\iOR u\pt\psi dyd\eta d\tau\\
=&\int_{s}^{\tau}\iOR(\eta\cdot\nabla_y\psi-h^2 y\cdot\nabla_{\eta}\psi)u dyd\eta d\tau
-\int_{s}^{\tau}\int_{\gamma}\psi u d\gamma d\tau\\
&+\int_{s}^{\tau}\iOR\mut\cos^2(h\tau)(Lu)\psi dyd\eta d\tau+\int_{s}^{\tau}\iOR\cos^2(h\tau)g\psi dyd\eta d\tau\\
:=&I_1+I_2+I_3+I_4,
\end{split}
\end{equation}
where $d\ga=(\eta\cdot n_y)dS_yd\eta$, and
$\psi\in C^{\infty}\cap L^2((0,\pi/2h)\times\O\times\R^3)$ is a test function.

We now focus on
the treatment of \eqref{eq: weak form for linear:1st}, which is divided into the following six steps.
In Step 1-Step 3, some useful estimates are derived for different choices of test function $\psi$. In Step 4-Step 6, 
based on Step 1-Step 3, the functions $a, b, c$
in \eqref{eq: u of abc} are dealt with.

\vskip 0.2 true cm

\textbf{Step 1. Choosing the test function $\psi=\vphi(y)\tmp\chi_0$ in
\eqref{eq: weak form for linear:1st}}

\vskip 0.2 true cm

 Direct computation yields
\begin{equation}\label{4-1}
\eta\cdot\nabla_y\psi-h^2y\cdot\nabla_{\eta}\psi
=\sum_{k=1}^{3}\p_{k}\vphi\tmp\chi_k.
\end{equation}
Let $[s,\tau]=[\tau,\tau+\ep]$ in \eqref{eq: weak form for linear:1st}.
Then it follows from \eqref{eq: weak form for linear:1st}, \eqref{4-1} and \eqref{eq: u of abc} that
the left-hand side (will be briefly written as {\bf LHS}) of equation (\ref{eq: weak form for linear:1st}) now becomes
\begin{equation}\label{4-2}
\text{\bf{LHS}}=\int_{\O}\left\{a(\tau+\ep)-a(\tau)\right\}\vphi dy.
\end{equation}
By (\ref{4-1}),
\begin{equation*}
I_1=\int_{\tau}^{\tau+\ep}\int_{\O}(b\cdot\nabla_y\vphi)  dy d\tau.
\end{equation*}
And $I_2=0$, $I_4=0$ by the boundary condition in (\ref{data for u in L2}) and
the assumption of $g$ in (\ref{mass conservation for g}), respectively. In addition,  the fact
of $<\chi_0,Lu>=0$ yields $I_3=0$. In this case, taking the difference quotient in \eqref{eq: weak form for linear:1st}
as $\ep\to0$ and using \eqref{4-2} yield
\begin{equation}\label{eq: partial a}
\int_{\O}\vphi\pt a dy=\int_{\O}(b\cdot\nabla_y\vphi) dy.
\end{equation}
Let $\vphi\equiv 1$ in \eqref{eq: partial a}, then
\begin{equation}\label{4-3}
\int_{\O}\pt ady=0.
\end{equation}
On the other hand, for $\vphi\in H^1(\O)$, we have that from \eqref{eq: partial a},
\begin{equation*}
|\int_{\O}\vphi\pt a dy|\lesssim\|b\|\|\vphi\|_{H^1}.
\end{equation*}
This leads to
\begin{equation}\label{4-4}
\|\pt a(\tau)\|_{H^{-1}_0}\lesssim\|b\|.
\end{equation}
For fixed  $\tau\in(0,\pi/2h)$,
by \eqref{4-3}-\eqref{4-4} and the standard elliptic theory, there exists a unique weak solution
$\Phi_a$ to the problem
\begin{equation*}
\begin{cases}
-\Delta\Phi_a   =\pt a(\tau),  & \text{in $\O$}, \\
\frac{\p\Phi_a}{\p n}  =0,     & \text{on $\p\O$}, \\
\int_{\O}\Phi_a   dy=0.        &
\end{cases}
\end{equation*}
Moreover,
\begin{equation}\label{estimate for pt a}
\|\Phi_a(\tau)\|_{H^1}\lesssim\|\pt a(\tau)\|_{H^{-1}_0}\lesssim\|b\|.
\end{equation}

\vskip 0.2 true cm

\textbf{Step 2. Choosing the test function $\psi_j=\vphi(y)\tmp\chi_j$ ($j=1,2,3$) in
\eqref{eq: weak form for linear:1st}}

\vskip 0.2 true cm

In this case, we have
\begin{equation}\label{4-6}
\begin{aligned}
&\eta\cdot\nabla_y\psi_j-h^2y\cdot\nabla_{\eta}\psi_j \\
=&\p_{j}\vphi\tmp(\chi_0+\frac{\sqrt{6}}{3}\chi_4)
-h^2y^j\vphi\tmp\chi_0+\sum_{k=1}^{3}\p_{k}\vphi\tmp B_{kj}.
\end{aligned}
\end{equation}
When $[s,\tau]=[\tau,\tau+\ep]$
is chosen in \eqref{eq: weak form for linear:1st}, the left-hand side of (\ref{eq: weak form for linear:1st}) is
\begin{equation}\label{4-7}
\text{\bf LHS}=\iO\left\{ b^j(\tau+\ep)- b^j(\tau)\right\}\vphi dy.
\end{equation}
Meanwhile, by \eqref{4-6},
\begin{equation}\label{4-8}
\begin{aligned}
I_1= \int_{\tau}^{\tau+\ep}\iO
\left\{\p_{j}\vphi(a+\frac{\sqrt{6}}{3}c-\frac{\sqrt{6}}{3}q)
-h^2y^j\vphi a\right\} dy d\tau
+\int_{\tau}^{\tau+\ep}\iO\sum_{k=1}^{3}\p_{k}\vphi<B_{kj},d>dy d\tau
\end{aligned}
\end{equation}
and $I_4=0$. In addition, from the fact that $<Lu,\chi_j>=0$, we have $I_3=0$.
Consequently, taking the difference quotient as $\ep\to 0$ in \eqref{eq: weak form for linear:1st}
and using \eqref{4-7}-\eqref{4-8} yield
\begin{equation}\label{4-9}
\begin{aligned}
\int_{\O}\vphi\pt b^j dy= & \int_{\O}
\left\{\p_{j}\vphi(a+\frac{\sqrt{6}}{3}c-\frac{\sqrt{6}}{3}q)
-h^2y^j\vphi a\right\}dy\\
& +\int_{\O}\sum_{k=1}^{3}\p_{k}\vphi<B_{kj},d>dy
+\int_{\gamma}u\tmp\vphi\chi_jd\gamma.
\end{aligned}
\end{equation}
For fixed $\tau>0$, let $\Phi_b=(\Phi_b^1,\Phi_b^2,\Phi_b^3)$ satisfy
\begin{equation}\label{eq: elliptic partial tau Phib}
\left\{ \enspace
\begin{aligned}
-\Delta\Phi_b  & =\pt b(\tau)      & \qquad  \text{in $\O$}, \\
\Phi_b\cdot n        & =0           & \qquad  \text{on $\p\O$}, \\
\p_n\Phi_b    & =(\Phi_b\cdot n)n   & \qquad  \text{on $\p\O$}.
\end{aligned}
\right.
\end{equation}
The existence of $\Phi_b$ is given in Appendix A.
Choosing $\vphi=\Phi_b^j$ in \eqref{4-9}, after summation of $j=1,2,3$, the last term on the right hand of (\ref{4-9}) is
\bes
  \sum_{j=1}^{3}\int_{\gamma}u\tmp\Phi_b^j\chi_jd\gamma
    =2\int_{\p\Omega}\int_{\eta\cdot n>0}(\Phi_b\cdot n)(\eta\cdot n)u\tmp d\gamma=0.
\ees
So
\begin{align*}
\sum_{j=1}^{3}\int_{\O}|\nabla_y\Phi_b^j|^2dy= & \sum_{j=1}^{3}\int_{\O}-\Delta\Phi_b^j\Phi_b^j dy
=\sum_{j=1}^{3}\int_{\O}\Phi_b^j\pt b^j dy\\
\lesssim & \|\nabla_y\Phi_b^j\|(\|a\|+h^2\|a\|+\|c\|)
+h^2\|\Phi_b^j\|\|a\|
+\|\nabla_y\Phi_b^j\|\|d\|.
\end{align*}
This, together with $\|\Phi_b^j\|\lesssim\|\nabla_y\Phi_b^j\|$, yields
\begin{equation}\label{estimate for pt b}
\|\Phi_b^j(\tau)\|_{H^1}\lesssim (\|a\|+\|c\|+\|d\|).
\end{equation}

\vskip 0.2 true cm

\textbf{Step 3. Choosing the test function $\psi=\vphi(y)\tmp\chi_4$  in
\eqref{eq: weak form for linear:1st}}

\vskip 0.2 true cm

In this case, one has
\begin{equation}\label{4-A}
\eta\cdot\nabla_y\psi-h^2 y\cdot\nabla_{\eta}\psi
=\sum_{k=1}^{3}\frac{\sqrt{6}}{3}(\p_{k}\vphi-h^2 y^k\vphi)\tmp\chi_k
+\sum_{k=1}^{3}\frac{\sqrt{15}}{3}\p_{k}\vphi\tmp A_k.
\end{equation}
Let $[s,\tau]=[\tau,\tau+\ep]$ in \eqref{eq: weak form for linear:1st}. Then
it follows from \eqref{eq: weak form for linear:1st} and \eqref{eq: u of abc}
that
\begin{equation*}
\text{\bf LHS}=\int_{\O}\{c(\tau+\ep)-q(\tau+\ep)-c(\tau)-q(\tau)\}\vphi dy.
\end{equation*}
Meanwhile, by \eqref{4-A},
\begin{equation*}
I_1=\int_{\tau}^{\tau+\ep}\int_{\O}\left\{\frac{\sqrt{6}}{3}(b\cdot\nabla_y\vphi-h^2 b\cdot y\vphi)
+\sum_{k=1}^{3}\frac{\sqrt{15}}{3}\p_{k}\vphi<A_k,d>\right\}dyd\tau
\end{equation*}
and $I_2=I_4=0$. In addition, the fact that $<Lu,\chi_4>=0$ derives $I_3=0$. Consequently, as in Step 1 and Step 2,
we have
\begin{equation}\label{4-B}
\int_{\O}\vphi\pt c dy= \int_{\O}\frac{\sqrt{6}}{3}b\cdot\nabla_y\vphi dy
+\int_{\O}\sum_{k=1}^{3}\frac{\sqrt{15}}{3}\p_{k}\vphi<A_k,d> dy.
\end{equation}
Choosing $\vphi\equiv 1$ in \eqref{4-B} yields
\begin{equation*}
\int_{\O}\pt c dy=0.
\end{equation*}
Thus, for fixed $\tau>0$, we can choose $\vphi=\Phi_c$ such that
\begin{equation}\label{4-C}
\begin{cases}
-\Delta\Phi_c          =\pt c(\tau)    &  \text{in $\O$}, \\
\frac{\p\Phi_c}{\p n}  =0              &  \text{on $\p\O$}, \\
\int_{\O}\Phi_c       dy=0.            &
\end{cases}
\end{equation}
This leads to
\begin{align*}
  \int_{\O}|\nabla_y\Phi_c(\tau)|^2dy= & \int_{\O}-\Delta\Phi_c\Phi_c dy=\int_{\O}\Phi_c\pt c dy\\
     \lesssim & \|\nabla_y\Phi_c\|(\|b\|+\|d\|).
\end{align*}
Thus, we have
\begin{equation}\label{estimate for Pt_c}
\|\Phi_c(\tau)\|_{H^1}\lesssim (\|b\|+\|d\|).
\end{equation}

\vskip 0.2 true cm

Before continuing to Step 4-6,
we rewrite (\ref{eq: weak form for linear:1st}) in the following form
\begin{equation}\label{eq: weak form_2nd}
\begin{aligned}
& -\int_{s}^{\tau}\iOR(\eta\cdot\nabla_y\psi-h^2 y\cdot\nabla_{\eta}\psi)u dy d\eta d\tau\\
= & \iOR\left\{-\psi u(\tau)+\psi u(s)\right\} dy d\eta
+\int_{s}^{\tau}\iOR u\pt\psi  dy d\eta d\tau
-\int_{s}^{\tau}\int_{\gamma}\psi u d\gamma d\tau\\
& +\int_{s}^{\tau}\iOR\mut\cos^2(h\tau)(Lu)\psi  dy d\eta d\tau
+\int_{s}^{\tau}\iOR\cos^2(h\tau)g\psi  dy d\eta d\tau\\
:= & \left\{G_{\psi}(\tau)-G_{\psi}(s)\right\}+J_1+J_2+J_3+J_4.
\end{aligned}
\end{equation}
We will consider the weak formulation (\ref{eq: weak form_2nd}) instead of (\ref{eq: weak form for linear:1st}) in the following.

\vskip 0.2 true cm

\textbf{Step 4. Estimates of $c$}

\vskip 0.2 true cm

For fixed $\tau>0$, by \eqref{eq: conserve c} let $\vphi_c$ be a solution of the following problem
\begin{equation}\label{eq: vphi c}
\begin{cases}
-\Lap_y                  \vphi_c=c(\tau), & \text{in $\O$}, \\
\frac{\p\vphi_c}{\p n}   =0,             & \text{on $\p\O$}, \\
\int_{\O}\vphi_c        dy=0.            &
\end{cases}
\end{equation}
Choosing
\bes\ba
\psi= & \psi_c=:\sum_{j=1}^{3}\p_j\vphi_c\eta^j(|\eta|^2-5)M^{1/2}\\
    = & \sum_{j=1}^{3}\sqrt{10}\p_j\vphi_c\tmp A_j,
\ea\ees
in (\ref{eq: weak form_2nd}).
In addition, direct computation  yields
\bes
\begin{aligned}
&\eta\cdot\nabla_y\psi-h^2y\cdot\nabla_{\eta}\psi \\
= & \sum_{j,k=1}^{3}\p^2_{jk}\vphi_c\eta^k\eta^j(|\eta|^2-5)M^{1/2}
-h^2\sum_{j,k=1}^{3}y^k\p_j\vphi_c\tmp
\left(2B_{jk}+\frac{5\sqrt{6}}{3}\delta_{jk}\chi_4\right) \\
=&\frac{5\sqrt{6}}{3}\Lap\vphi_c\chi_4
+\sum_{j,k=1}^{3}\p^2_{jk}\vphi_c\tmp
(I-P)(\eta^k\eta^j(|\eta|^2-5)\mu^{1/2}) \\
&-h^2\sum_{j,k=1}^{3}y^k\p_j\vphi_c\tmp
\left(2B_{jk}+\frac{5\sqrt{6}}{3}\delta_{jk}\chi_4\right).
\end{aligned}
\ees
Thus, the left-hand side of (\ref{eq: weak form_2nd}) is
\begin{equation}\label{4-a}
\begin{aligned}
\text{\bf LHS}= & -\frac{5\sqrt{6}}{3}\int_{s}^{\tau}\iO(c-q)c dy d\tau
+\int_{s}^{\tau}\iO\p^2_{jk}\vphi_c
<d,\eta^k\eta^j(|\eta|^2-5)\mu^{1/2}> dy d\tau \\
& -h^2\sum_{j,k=1}^{3}\int_{s}^{\tau}\iO y^k\p_j\vphi_c
\left(\frac{5\sqrt{6}}{3}\delta_{jk}c+2<d,B_{jk}>\right) dy d\tau \\
:= & -\frac{5\sqrt{6}}{3}\int_{s}^{\tau}\iO(c-q)c dy d\tau+E_1,
\end{aligned}
\end{equation}
where, for any $\ve>0$,
\begin{equation}\label{4-b}
|E_1| \le (\ve^2+h^2)\|c\|^2+\frac{1}{\ep^2}\|d\|^2,
\end{equation}
and
\begin{equation}\label{4-c}
\begin{aligned}
|\int_{s}^{\tau}\iO cq dy d\tau|
\leq \|c\|\|q\|\leq h^2\|c\|^2+h^2\|a\|^2.
\end{aligned}
\end{equation}
On the other hand, by the fact that $\Phi_c=\pt\vphi_c$ and estimate of $\Phi_c$ in (\ref{estimate for Pt_c}),
we have that for any $\ve>0$,
\begin{equation}\label{4-d}
\begin{aligned}
|J_1| & \lesssim \sum_{j=1}^{3}\ist\iO|\p^2_{\tau j}\vphi_c<A_j,d>| dy d\tau\\
& \lesssim \int_{s}^{\tau}\|\pt\vphi_c\|_{H^1}\|d\| d\tau\\
& \lesssim \ist(\|b\|+\|d\|)\|d\| d\tau\\
& \lesssim \ep^2\|b\|^2+\frac{1}{\ep^2}\|d\|^2.
\end{aligned}
\end{equation}
By the boundary conditions of $\vphi_c$ and $u$, one has that
\bes
\begin{aligned}
& \int_{\p\O\times\R^3}\psi ud\gamma \\
= & \sum_{j=1}^{3}\int_{\p\O}\p_j\vphi_c
\left(\int_{\eta\cdot n_y>0}+\int_{\eta\cdot n_y<0}\right)
(\eta\cdot n_y)\eta^j(|\eta|^2-5)\mu^{1/2}u(\eta) d\gamma \\
= & \sum_{j=1}^{3}\int_{\p\O}\p_j\vphi_c\int_{\eta\cdot n_y>0}
(\eta\cdot n_y)\eta^j(|\eta|^2-5)\mu^{1/2}u(\eta) d\gamma  \\
& +\sum_{j=1}^{3}\int_{\p\O}\p_j\vphi_c
\int_{\eta\cdot n_y>0}(-\eta\cdot n_y)
(\eta^j-2(\eta\cdot n_y)n_y^j)
(|\eta|^2-5)\mu^{1/2}u(\eta) d\gamma  \\
= & 2\int_{\p\O}\frac{\p\vphi_c}{\p n_y}
\int_{\eta\cdot n_y>0}
(\eta\cdot n_y)^2(|\eta|^2-5)\mu^{1/2}u  d\gamma  \\
= & 0.
\end{aligned}
\ees
Therefore,
\begin{equation}\label{4-e}
\begin{aligned}
J_2=0.
\end{aligned}
\end{equation}
In addition, we obtain that for any $\ve>0$,
\begin{align}\label{4-f}
|J_3|\lesssim & \ist\iO\cos^2(h\tau)|(Ld)\psi| dyd\tau \no\\
\lesssim & \ist\cos^2(h\tau)\|d\|\|\nabla_y\vphi_c\| d\tau\no\\
\lesssim & \ep\|c\|^2+\epn\|d\|^2,\no\\
|J_4|\lesssim & \ist\cos^2(h\tau)\|g\|\|\psi\|  dyd\tau \no\\
\lesssim & \ist\cos^2(h\tau)\|g\|\|c\| d\tau\no\\
\lesssim & \ep^2\|c\|^2+\frac{1}{\ep^2}\|g\|^2.
\end{align}
By choosing $\ep>0$ and $h>0$ small, it follows from \eqref{4-a}-\eqref{4-f} that
\begin{equation}\label{4-g}
\begin{aligned}
\|c\|^2\lesssim  (G_{c}(\tau)-G_c(s)) +\ep^2\|a\|^2+\ep^2\|b\|^2
+\frac{1}{\ep^2}\|d\|^2+\frac{1}{\ep^2}\|g\|^2.
\end{aligned}
\end{equation}

\vskip 0.2 true cm

\textbf{Step 5. Estimate of $b$}

\vskip 0.2 true cm

For fixed $\tau>0$, let $\vphi_b=(\vphi_b^1,\vphi_b^2,\vphi_b^3)$
be the solution of the following problem:
\begin{equation}\label{eq: vphi b}
\begin{cases}
-\Lap_y  \vphi_b=  b(\tau) &  \text{ in $\O$}, \\
\vphi_b\cdot n =  0        &  \text{ on $\p\O$}, \\
\p_n\vphi_b=  (\p_n\vphi_b\cdot n)n &  \text{ on $\p\O$},
\end{cases}
\end{equation}
where the existence of solution $\vphi_b$ is proved in Appendix A.
Set
\bes\ba
\psi= & \psi_b
=\ds\sum_{i,j=1}^{3}\p_j\vphi_b^i\eta^i\eta^jM^{1/2}
-\sum_{i=1}^{3}\p_i\vphi_b^i\frac{|\eta|^2-1}{2}M^{1/2} \\
= & \sum_{i,j=1}^{3}\p_j\vphi_b^iB_{ij}
-\sum_{i=1}^{3}\frac{\sqrt{6}}{6}\p_i\vphi_b^i\chi_4,
\ea\ees
in (\ref{eq: weak form_2nd}), we then have that by direct computation,
\begin{equation*}
\begin{aligned}
& \eta\cdot\nabla_y\psi-h^2y\cdot\nabla_{\eta}\psi \\
= & \sum_{i=1}^{3}\Lap\vphi_b^i\tmp\chi_i
+\sum_{i,j,k=1}^{3}\p^2_{jk}\vphi_b^i\tmp
(I-P)(\eta^i\eta^j\eta^k\mu^{1/2}) \\
&-h^2\sum_{i,j,k=1}^{3}y^k\p_j\vphi_b^i\tmp
(\delta_{ik}\chi_j+\delta_{jk}\chi_i)
+h^2\sum_{i,k=1}^{3}y^k\p_i\vphi_b^i\tmp\chi_k.
\end{aligned}
\end{equation*}
Thus, the left-hand side of (\ref{eq: weak form_2nd}) becomes
\begin{equation}\label{0-0}
\begin{aligned}
\text{\bf LHS}= & \sum_{i=1}^{3}\ist\iO b^i\Lap\vphi_b^i dyd\tau
+\sum_{i,j,k=1}^{3}\ist\iO\p^2_{jk}\vphi_b^i
<\eta^i\eta^j\eta^k\mu^{1/2},d>  dyd\tau\\
& -h^2\sum_{i,j,k=1}^{3}\ist\iO \bigl(y^k\p_j\vphi_b^i
(\delta_{ik}b^j+\delta_{jk}b^i)
+h^2\sum_{i,k=1}^{3}y^k\p_i\vphi_b^ib^k\bigr) dyd\tau \\
:= & \sum_{i=1}^{3}\ist\iO b^i\Lap\vphi_b^i dyd\tau+E_2,
\end{aligned}
\end{equation}
where, for any $\ve>0$,
\begin{equation}\label{0-1}
|E_2| \le (\ep^2+h^2)\|b\|^2+\frac{1}{\ep^2}\|d\|^2.
\end{equation}
Next we estimate $J_i$ $(1\le i\le 4)$ in \eqref{eq: weak form_2nd}.
From the fact that $\Phi_b=\pt\vphi_b$ and the estimate of $\Phi_b$ in (\ref{estimate for pt b}), we have
that for any $\ve>0$,
\begin{equation}\label{0-2}
\begin{aligned}
|J_1|\lesssim & \ist(\|c\|+\|d\|)\|\pt\nabla_y\vphi_b\| d\tau\\
\lesssim & \ep^2\|a\|^2
+\frac{1}{\ep^2}\|c\|^2+\frac{1}{\ep^2}\|d\|^2+\epn\|g\|^2.
\end{aligned}
\end{equation}
In addition,
\begin{equation}\label{R-0}
\begin{aligned}
\int_{\gamma}\psi ud\gamma= & \int_{\gamma}(\eta\cdot n_y)
\left(\sum_{i,j=1}^{3}\p_j\vphi_b^i\eta^i\eta^jM^{1/2}
-\p_i\vphi_b^i\frac{|\eta|^2-5}{2}M^{1/2}\right)u d\gamma\\
= & \int_{\gamma}(\eta\cdot n_y)
\sum_{i,j=1}^{3}\p_j\vphi_b^i\eta^i\eta^jM^{1/2}u d\gamma.
\end{aligned}
\end{equation}
For fixed $y\in\p\O$, set
\bes
J_{20}:=\sum_{i,j=1}^{3}\int_{\R^3}(\eta\cdot n_y)
\p_j\vphi_b^i\eta^i\eta^jM^{1/2}u d\eta.
\ees
Through a coordinate rotation, we may assume that $n_y=(1,0,0)$. Then
\bes
\begin{aligned}
J_{20}= & \p_1\vphi_b^1\int_{\R^3}(\eta^1)^3M^{1/2}ud\eta
+\sum_{i=2}^{3}\p_1\vphi_b^i\int_{\R^3}(\eta^1)^2\eta^iM^{1/2}ud\eta\\
&+\sum_{j=2}^{3}\p_j\vphi_b^1\int_{\R^3}(\eta^1)^2\eta^jM^{1/2}ud\eta
+\sum_{i,j=2}^{3}\p_j\vphi_b^i\int_{\R^3}\eta^1\eta^i\eta^jM^{1/2}ud\eta\\
& :=J_{21}+J_{22}+J_{23}+J_{24}.
\end{aligned}
\ees
From the boundary condition of $u$ and $n_y=(1,0,0)$, we have $u(\eta^1,\eta^2,\eta^3)=u(-\eta^1,\eta^2,\eta^3)$.
So $J_{21}=J_{24}=0$. On the other hand, since the second boundary condition of $\vphi_b$ in (\ref{eq: vphi b}) gives $\p_1\vphi_b^i=0$ for $i\neq 1$,
$J_{22}=0$ holds;  since the first boundary condition of $\vphi_b$ in (\ref{eq: vphi b}) gives $\vphi_b^1=0$ on $\p\O$,
$\p_j\vphi_b^1=0$ for $j\neq 1$ and further $J_{23}=0$ hold. Thus by \eqref{R-0} and the expression of $J_2$ we obtain
\begin{equation}\label{0-3}
\begin{aligned}
J_2=0.
\end{aligned}
\end{equation}
Meanwhile, it follows from direct computation that
\begin{equation}\label{0-4}
\begin{aligned}
|J_3|\lesssim & \ist \|d\|\|b\|d\tau\\
\lesssim & \ep^2\|b\|^2+\frac{1}{\ep^2}\|d\|^2,
\end{aligned}
\end{equation}
and
\begin{equation}\label{0-5}
\begin{aligned}
|J_4|\lesssim & \ist\|g\|\|b\|d\tau \\
\lesssim & \ep^2\|b\|^2+\frac{1}{\ep^2}\|g\|^2.
\end{aligned}
\end{equation}
Substituting \eqref{0-0}-\eqref{0-5} into \eqref{eq: weak form_2nd} yields that for small $\ve>0$,
\begin{equation}\label{0-6}
\begin{aligned}
\|b\|^2\lesssim  (G_b(\tau)-G_b(s))+\ep^2\|a\|^2+\frac{1}{\ep^2}\|c\|^2+\frac{1}{\ep^2}\|d\|^2+\frac{1}{\ep^2}\|g\|^2.
\end{aligned}
\end{equation}

\vskip 0.2 true cm

\textbf{Step 6. Estimates of $a$}

\vskip 0.2 true cm

For fixed $\tau>0$, by \eqref{eq: conserve a} let $\vphi_a$ be the solution of the following problem:
\begin{equation}\label{eq: vphi a}
\begin{cases}
-\Lap_y                  \vphi_a=a(\tau) & \text{in $\O$}, \\
\frac{\p\vphi_a}{\p n}   =0              & \text{on $\p\O$}, \\
\int_{\O}\vphi_a        dy=0.            &
\end{cases}
\end{equation}
Choosing
\bes\ba
\psi= & \psi_a=\sum_{j=1}^{3}\p_j\vphi_a\eta^j(|\eta|^2-10)M^{1/2} \\
= & \sum_{j=1}^{3}\p_j\vphi_a\tmp(\sqrt{10}A_j-5\chi_j).
\ea\ees
 Then direct computation yields
\begin{equation*}
\begin{aligned}
& \eta\cdot\nabla_y\psi-h^2y\cdot\nabla_{\eta}\psi \\
=&\sum_{j,k=1}^{3}\p^2_{jk}\vphi_a\eta^k\eta^j(|\eta|^2-10)M^{1/2}\\
&-h^2\sum_{j,k=1}^{3}y^k\p_j\vphi_a
\{\delta_{jk}(|\eta|^2-10)+2\eta^j\eta^k\}M^{1/2} \\
= & -5\Lap\vphi_a\tmp\chi_0
+\sum_{j,k=1}^{3}\p^2_{jk}\vphi_a\tmp
(I-P)\{\eta^k\eta^j(|\eta|^2-10)\mu^{1/2}\}\\
& -h^2\sum_{j,k=1}^{3}y^k\p_j\vphi_a\tmp
\left(\delta_{jk}(\frac{5\sqrt{6}}{3}\chi_4-5\chi_0)+2B_{jk}\right).
\end{aligned}
\end{equation*}
Thus, the left-hand side of (\ref{eq: weak form_2nd}) becomes
\begin{equation}\label{a-0}
\begin{aligned}
\text{\bf LHS}=&5\|a\|^2+\sum_{j,k=1}^{3}\ist\int_{\O}\p^2_{jk}\vphi_a
<\eta^k\eta^j(|\eta|^2-10)\mu^{1/2},d> dyd\tau\\
& -h^2\sum_{j,k=1}^{3}\ist\iO y^k\p_j\vphi_a
\left(\frac{5\sqrt{6}}{3}\delta_{jk}c-5\delta_{jk}a+2<B_{jk},d>\right)dy d\tau \\
:= & 5\|a\|^2+E_3,
\end{aligned}
\end{equation}
where
\begin{equation}\label{a-1}
\begin{aligned}
|E_3|\le (\ep^2+h^2)\|a\|^2+h^2\|c\|^2+\frac{1}{\ep^2}\|d\|^2.
\end{aligned}
\end{equation}
Next we estimate $J_i$ $(1\le i\le 4)$ in \eqref{eq: weak form_2nd}.
From the fact that $\Phi_a=\pt\vphi_a$ and the estimate of $\Phi_a$ in \eqref{estimate for pt a}, we have
\begin{equation*}
\begin{aligned}
|J_1|\lesssim & \ist(\|b\|+\|d\|)\|\pt\vphi_a\|_{H^1} d\tau\\
\lesssim & \ist(\|b\|+\|d\|)\|b\| d\tau\\
\lesssim & \|b\|^2+\|d\|^2.
\end{aligned}
\end{equation*}
Similar to the treatment in Step 4, we have
\begin{equation*}
\int_{\p\O\times\R^3}\psi ud\gamma=
2\int_{\p\O}\frac{\p\vphi_c}{\p n_y}\int_{\eta\cdot n_y>0}(\eta\cdot n)^2(|\eta|^2-10)\mu^{1/2}d =0
\end{equation*}
and thus
\begin{equation}\label{a-2}
\begin{aligned}
J_2=0.
\end{aligned}
\end{equation}
In addition, we have
\begin{equation}\label{a-3}
\begin{aligned}
|J_3|\lesssim & \ist\|Ld\|\|\psi\|d\tau\\
\lesssim & \ep^2\|a\|^2+\frac{1}{\ep^2}\|d\|^2,
\end{aligned}
\end{equation}
and
\begin{equation}\label{a-4}
\begin{aligned}
|J_4|\lesssim & \ist\|g\|\|\psi\| d\tau\\
\lesssim & \ep^2\|a\|^2+\frac{1}{\ep^2}\|g\|^2.
\end{aligned}
\end{equation}
It follows from \eqref{a-0}-\eqref{a-4} that
\begin{equation}\label{a-5}
\begin{aligned}
\|a\|^2\lesssim  (G_a(\tau)-G_a(s))+\ep^2\|c\|^2+\frac{1}{\ep^2}\|b\|^2+\frac{1}{\ep^2}\|d\|^2+\frac{1}{\ep^2}\|g\|^2.
\end{aligned}
\end{equation}
Using the fact that $\|d\|\leq\|d\|_{\nu}$ and the asymptotic behavior of $\nu$ in (\ref{bound for nu}), we have
\begin{equation}\label{ineq: norm nu with norm}
\|a\tmp\|_{\nu}=\|a\nu^{1/2}\tmp\|\lesssim \|a\|.
\end{equation}
Combining (\ref{a-5}), (\ref{0-6}) and (\ref{4-g}) with (\ref{ineq: norm nu with norm}) yields
\begin{equation}\label{estimates: Pu in terms of G and d}
\ist\|Pu\|_{\nu}^2d\tau\lesssim (G(\tau)-G(s))+\ist\|(I-P)u\|_{\nu}^2d\tau+\|g\|^2,
\end{equation}
where
\be\label{estimate: G by u}
\begin{aligned}
|G(\tau)|= & |G_a+G_b+G_c| \\
= & |\iO-(\psi_a+\psi_b+\psi_c)u(\tau)dy| \\
\lesssim & \|u(\tau)\|^2.
\end{aligned}
\ee

Next we derive the $L^2$-decay of solution $u$ to the linear Boltzmann equation \eqref{eq: linear u with g}.
For some constant $\la>0$ to be determined
later on, set
\bes
U(\tau)=u(\tau)e^{\la\al(\tau)}.
\ees
Then $U$ satisfies
\begin{equation}\label{R-1}
\Lat U
=\mut\cos^2(h\tau)LU+\la\cos^2(h\tau)U+\cos^2(h\tau)g e^{\la\al(\tau)}.
\end{equation}
Multiplying both sides of \eqref{R-1} by $U$ and integrating with respect to the variable $(\tau,y,\eta)$,
and using the dissipation property of $L$ in (\ref{negetivity of L}), we arrive at
\begin{equation}\label{eq: U-U}
\begin{aligned}
& \|U(\tau)\|^2-\|U(s)\|^2 \\
= & 2\ist\iO\mut\cos^2(h\tau)<LU,U>dyd\tau
+2\la\ist\cos^2(h\tau)\|U\|^2d\tau\\
& +2\ist\iO\cos^2(h\tau)<g e^{\la\al(\tau},U>dyd\tau \\
\lesssim & -\ist\cos^2(h\tau)\|(I-P)U\|_{\nu}^2d\tau
+\la\ist\cos^2(h\tau)\|U\|^2d\tau \\
& +\frac{1}{\la}\ist\cos^2(h\tau) e^{2\la\al(\tau)}\|g\|^2d\tau.
\end{aligned}
\end{equation}
For $n\geq0$, set
\bes
\begin{aligned}
\tau_n:= & \frac{\arctan(hn)}{h}\in[0,\pi/2h), \\
m_n:= & \min_{\tn\leq s\leq\tnp}
\{\cos^2(hs)e^{2\la\al(s)}t\},\\
M_n:= & \max_{\tn\leq s\leq\tnp}
\{\cos^2(hs)e^{2\la\al(s)}\}.
\end{aligned}
\ees
Let $s=\tn$ and $\tau=\tnp$ in (\ref{eq: U-U}). Then we have
\be\label{eq: U tau s}
\begin{aligned}
& \|U(\tnp)\|^2-\|U(\tn)\|^2 \\
\lesssim & -m_n\int_{\tn}^{\tnp}\|(I-P)u\|_{\nu}^2d\tau
+\la M_n\int_{\tn}^{\tnp}\|u\|^2d\tau
 +\frac{M_n}{\la}\int_{\tn}^{\tnp}e^{\la\al(\tau)}\|g\|^2d\tau.
 \end{aligned}
\ee
Furthermore, (\ref{eq: U tau s})$+\{$(\ref{estimates: Pu in terms of G and d})
$\times\ep m_n$\} for some small $\ep>0$ gives
\bes
\begin{aligned}
& \{\|U(\tnp)\|^2-\ep m_n G(\tnp)\}
-\{\|U(\tn)\|^2-\ep m_n G(\tn)\} \\
\lesssim & -(m_n-\ep m_n-\la M_n)
\int_{\tn}^{\tnp}\|(I-P)u\|_{\nu}^2d\tau
-(\ep m_n-\la M_n)\int_{\tn}^{\tnp}\|Pu\|^2d\tau  \\
& +(\frac{M_n}{\la}+\ep m_n)
\int_{\tn}^{\tnp}e^{\la\al(\tau)}\|g\|^2d\tau.
\end{aligned}
\ees
Using the fact that $m_n\leq M_n\leq 4\e^{4\la} m_n\leq 4\e^{4}m_n$, we can choose $\ep\lesssim1/2$ and $\la\lesssim 1/16\e^4$ in the above.
In addition, by the estimate of $G$ in (\ref{estimate: G by u}), we then get
\bes
\|U(\tnp)\|^2-\|U(\tn)\|^2 \lesssim
\int_{\tn}^{\tnp}\cos^2(hs)e^{2\la\al(s)}\|g(s)\|^2ds.
\ees
Summing $n$ over $0,1,2...$ in the above, we conclude that
\begin{equation*}
\|u(\tau)\|^2\lesssim e^{-2\la\al(\tau)}\left\{\|u(0)\|^2
+\int_{0}^{\tau}e^{2\la\al(s)}\|g(s)\|^2ds\right\}.
\end{equation*}
Thus, the proof of Proposition 4.1 is completed. \qquad
\qquad \qquad \qquad \qquad \qquad \qquad \qquad \qquad \qquad $\square$

\section{$L^{\infty}$ decay of solutions to linear weighted Boltzmann equations}

Let $\phi$ be the weight function:
\bes
\phi=\phi_{\beta}(y,\eta)=(1+|\eta|^2+h^2|y|^2)^{\beta/2},
\ees
where $\beta>3/2$ is large enough.

Consider the equation of $w=\phi u$:
\begin{equation}\label{eq: w}
\Lat w+\mut(y)\cos^2(h\tau)\nu(\eta)w
=\mut(y)\cos^2(h\tau)K_{\phi}w+\cos^2(h\tau)\phi g,
\end{equation}
where
\bes
K_{\phi}f=\phi K\left(\frac{f}{\phi}\right).
\ees
The initial-boundary data of $w$ are
\begin{equation}\label{data for w in L infty}
\begin{cases}
 w(0,y,\eta)=\phi u_0(y,\eta),   &
 \text{$y\in \O,\,\eta\in\R^3$,}\\
 w(\tau,y,\eta)=w(\tau,y,\eta-2(\eta\cdot n_y)n_y),&
 \text{$y\in\p\O$, $\eta\cdot n_y<0$.}
\end{cases}
\end{equation}
Recall that $\al(\tau)=\int_{0}^{\tau}\cos^2(hs)ds$ and $\nu_0$ is defined in (\ref{bound for nu}),
we will prove the following conclusion.

\vskip 0.2 true cm

{\bf Proposition 5.1.} {\it
Let $w=w(\tau,y,\eta)$ be the mild solution of the linear Boltzmann
equation (\ref{eq: w}) with (\ref{data for w in L infty}). Then for some constant $\la\in(0,\nu_0/2)$, we have
\be\label{estimates: w L infty}
\|w(\tau)\|_{\infty}\lesssim
 e^{-\la\al(\tau)}\left(\|w_0\|_{\infty}
+\sup_{0<s<\tau}
\{e^{\la\al(s)}\|\nu^{-1}\phi g(s)\|_{\infty}\}\right).
\ee
}

To prove Proposition 5.1, we will express the solution $w(\tau)$ of  (\ref{eq: w}) by an integral along the backward trajectory.
After splitting the integration into several parts, we can bound the main part by the $L^2$
norm of $u(\tau)$ established in Section 4, and bound the remaining parts by $L^{\infty}$ norm of $w$.
Such a kind of estimate is called the $L^2-L^{\infty}$ estimate in \cite{MR2679358}.

For any fixed $(\tau,y,\eta)\in(0,\pi/2h)\times\O\times\R^3$,
set $(\tau_0',y_0',\eta_0')=(\tau,y,\eta)$. For $\tau_1\in(0,\tau)$,
by the formula (\ref{eq: y eta backward}), the backward trajectory $(y_1(\tau_1), \eta_1(\tau_1))$ is
\be\label{eq: Y1 and H1}
\begin{aligned}
y_1(\tau_1):=\sum_{k=0}^{\infty}
 1_{(\tau_{k+1}',\tau_k']}(\tau_1)Y(\tau_1;\tau_k',y_k',\eta_k'),\\
\eta_1(\tau_1):=\sum_{k=0}^{\infty}
 1_{(\tau_{k+1}',\tau_k']}(\tau_1)H(\tau_1;\tau_k',y_k',\eta_k').\\
\end{aligned}
\ee
To simplify the notation, along the trajectory $(y_1(\tau_1),\eta_1(\tau_1))$, we set
\bes
\begin{aligned}
\hal(\tau_1)= & \cos^2(h\tau_1), \\
\tal(\tau_1)= & \mut(y_1)\cos^2(h\tau_1), \\
\bal(\tau_1)= & \mut(y_1)\cos^2(h\tau_1)\nu(\eta_1).
\end{aligned}
\ees
By integrating along the trajectory, we have
\be\label{eq: w(tau)}
\begin{aligned}
w(\tau)= & w|_{\tau_1=0} e^{-\int_{0}^{\tau}\bal(\tau_1)d\tau_1}
+\int_{0}^{\tau}\tal(\tau_1) e^{-\int_{\tau_1}^{\tau}\bal(s_1)ds_1}
(K_{\phi}w)(\tau_1)d\tau_1\\
&+\int_{0}^{\tau}\hal(\tau_1) e^{-\int_{\tau_1}^{\tau}\bal(s_1)ds_1}
\phi g(\tau_1)d\tau_1 \\
:= & I_1+I_2+I_3.
\end{aligned}
\ee
Recall that in (\ref{def: set A}), for some $N>0$ large and $\kappa>0$ small, we have defined the set
\bes
\AehN:=\{(y_0,\eta_0)\in\O\times\R^3:\,
e_0-m_0\geq h^2+\ka^2,\, 2h\leq|\eta_0|\leq2N\}.
\ees
Next we start to estimate $w(\tau)$ in \eqref{eq: w(tau)}.

\vskip 0.2 true cm

\textbf{Step 1. Estimate of the main part $\|w(\tau)1_{\AehN}\|_{\infty}$}

\vskip 0.2 true cm

\begin{figure}[htbp]
\centering\includegraphics[width=10cm,height=6cm]{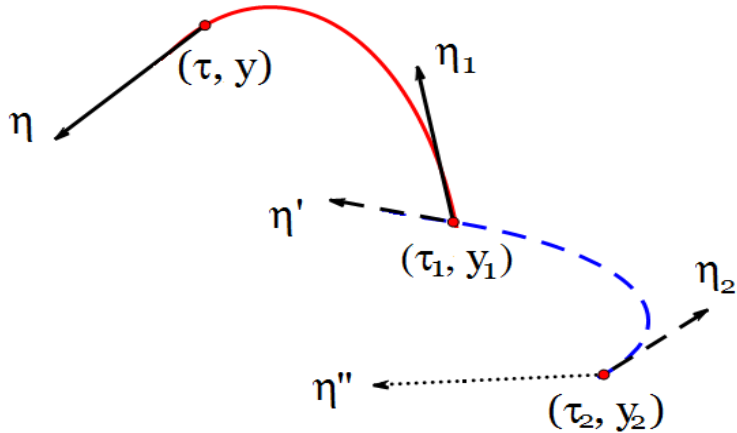}
\caption{\bf Trajectories in double integration}\label{fig:4}
\end{figure}

For any $(\tau,y,\eta)\in(0,\pi/2h)\times\AehN$, by $2\la\leq\nu_0\leq\nu(\eta_1)$ in (\ref{bound for nu}), we can bound $I_1$ in (\ref{eq: w(tau)}) by
\be\label{estimate: I1}
|I_1|\leq e^{-\la\al(\tau)}\|w_0\|_{\infty}.
\ee
For $I_3$, by the fact that
\bes
\int_{0}^{\tau}\nu(\eta_1)\hal(\tau_1)
 e^{-2\la\int_{0}^{\tau}\tal(s_1)ds_1}
\les  e^{-\la\al(\tau)},
\ees
we have
\be\label{estimate: I3}
|I_3|\les e^{-\la\al(\tau)}\sup_{0<s<\tau}
\{ e^{\la\al(s)}\|\nu^{-1}\phi g(s)\|_{\infty}\}.
\ee
Next we focus on the estimate of $I_2$. Note that
\bes
\begin{aligned}
(K_{\phi}w)(\tau_1)= & (K_{\phi}w)(\tau_1,y_1,\eta_1)
= \int_{\R^3}\kp(\eta_1,\eta')w(\tau_1,y_1,\eta')d\eta'.
\end{aligned}
\ees
Set $(\tau_0'',y_0'',\eta_0'')=(\tau_1,y_1,\eta')$, along the following backward trajectory
$(y_2(\tau_2), \eta_2(\tau_2))$
\be\label{eq: Y2 and H2}
\begin{aligned}
y_2(\tau_2):=\sum_{j=0}^{\infty}1_{(\tau_{j+1}'',\tau_j'']}(\tau_2)
Y(\tau_2;\tau_j'',y_j'',\eta_j''),\\
\eta_2(\tau_2):=\sum_{j=0}^{\infty}1_{(\tau_{j+1}'',\tau_j'']}(\tau_2)
H(\tau_2;\tau_j'',y_j'',\eta_j''),\\
\end{aligned}
\ee
$w(\tau,y_1,\eta')$ can be expressed as
\be\label{eq: w(tau1)}
\begin{aligned}
w(\tau_1)= & w|_{\tau_2=0} e^{-\int_{0}^{\tau_1}\bal(\tau_2)d\tau_2}
+\int_{0}^{\tau_1}\tal(\tau_2) e^{-\int_{\tau_2}^{\tau_1}\bal(s_2)ds_2}
(K_{\phi}w)(\tau_2)d\tau_2\\
&+\int_{0}^{\tau_1}\hal(\tau_2) e^{-\int_{\tau_2}^{\tau_1}\bal(s_2)ds_2}
\phi g(\tau_2)d\tau_2,
\end{aligned}
\ee
here we have used the simplified notations:
\bes
\begin{aligned}
\hal(\tau_2)= & \cos^2(h\tau_2), \\
\tal(\tau_2)= & \mut(y_2)\cos^2(h\tau_2), \\
\bal(\tau_2)= & \mut(y_2)\cos^2(h\tau_2)\nu(\eta_2).
\end{aligned}
\ees
One can see Figure \ref{fig:4} for these two trajectories in double integration
in \eqref{eq: w(tau1)}.

Substituting (\ref{eq: w(tau1)}) into the expression of $I_2$ in (\ref{eq: w(tau)}), we have
$I_2=I_{21}+I_{22}+I_{23}$, just
as in (\ref{eq: w(tau)}).
It follows from Lemma 3.3 that $K_{\phi}$ is a bounded operator, and then the following estimates of $I_{21}$ and $I_{23}$ hold (similar to $I_1$
in \eqref{estimate: I1} and $I_3$ in \eqref{estimate: I3}),
\be\label{estimate: I21}
|I_{21}|\les  e^{-\la\al(\tau)}\|w_0\|_{\infty},
\ee
and
\be\label{estimate: I23}
|I_{23}|\les  e^{-\la\al(\tau)}\sup_{0<s<\tau}
\{ e^{\la\al(s)}\|\nu^{-1}\phi g(s)\|_{\infty}\}.
\ee
The remaining part $I_{22}$ is
\be\label{eq: I22}
\begin{aligned}
I_{22}= & \int_{0}^{\tau}\tal(\tau_1) e^{-\int_{\tau_1}^{\tau}\bal(s_1)ds_1}
\int_{0}^{\tau_1}\tal(\tau_2) e^{-\int_{0}^{\tau_1}\bal(s_2)ds_2}\\
& \times \int_{\R^3\times\R^3}\kp(\eta_1,\eta')\kp(\eta_2,\eta'')
\Big\{ e^{\int_{0}^{\tau_2}\bal(s_2)ds_2}
w(\tau_2,y_2,\eta'')\Big\}d\eta'\eta''d\tau_1d\tau_2 \\
:= & I_{221}+I_{222}+I_{223},
\end{aligned}
\ee
where we have split the time-velocity integration into three cases, i.e.,
$I_{221}$ contains the integral domain: $|\eta|\geq N$ for some $N>0$ large enough;
$I_{222}$ contains the integral domain: $|\eta|\leq N$, $|\eta'|\geq 2N$ or $|\eta'|\leq 2N$,
$|\eta''|\geq 3N$; $I_{223}$  contains the integral domain: $|\eta|\leq N$, $|\eta'|\leq 2N$ and $|\eta''|\leq 3N$.
We next treat $I_{221}$, $I_{222}$ and $I_{223}$, respectively.

\vskip 0.2 true cm

\textbf{(a) The estimate of $I_{221}$}

\vskip 0.2 true cm

In this case, from the conservation of energy in (\ref{eq: e0}), one has
\be\label{estimate: eta1 by eta}
  |\eta_1|=(|\eta|^2+h^2|y|^2-h^2|y_1|^2)^{1/2}\geq N-1.
\ee
By Lemma 3.2 and (\ref{estimate: k}), we have that for $N$ large enough,
\bes
\int_{\R^3\times\R^3}\kp(\eta_1,\eta')\kp(\eta_2,\eta'')
d\eta'd\eta''\les\frac{1}{N}.
\ees
The double time integration can be estimated as follows
\bes
\begin{aligned}
& \int_{0}^{\tau}\tal(\tau_1) e^{-\int_{\tau_1}^{\tau}\bal(s_1)ds_1}
\int_{0}^{\tau_1}\tal(\tau_2) e^{-\int_{0}^{\tau_1}\bal(s_2)ds_2}
d\tau_1d\tau_2\\
\leq & \int_{0}^{\tau}\hal(\tau_1) e^{-\nu_0\al(\tau)+\nu_0\al(\tau_1)}
\int_{0}^{\tau_1}\hal(\tau_2) e^{-\nu_0\al(\tau_1)}d\tau_1d\tau_2\\
= & \frac{1}{2}\left(\al(\tau)\right)^2
 e^{-\nu_0\al(\tau)}\\
\les &  e^{-\la\al(\tau)},
\end{aligned}
\ees
which derives
\be\label{estimate: I221}
|I_{221}|\leq\frac{1}{N} e^{-\la\al(\tau)}
\sup_{0<s<\tau}\{ e^{\la\al(s)}|w(s)|\}.
\ee

\vskip 0.2 true cm

\textbf{(b) The estimate of $I_{222}$}

\vskip 0.2 true cm

In this case, similar to (\ref{estimate: eta1 by eta}), we have $|\eta_1-\eta'|\geq N-1$ or $|\eta''-\eta_2|\geq N-1$.
Then either of the following holds
\bes
|\kp(\eta_1,\eta')|\les \frac{1}{N}e^{-|\eta_1-\eta'|^2/8},
\ees
or
\bes
|\kp(\eta_2,\eta'')|\les \frac{1}{N}e^{-|\eta_2-\eta''|^2/8}.
\ees
From Lemma 3.2 and (\ref{estimate: k}), we arrive at
\bes
\int_{\R^3\times\R^3}\kp(\eta_1,\eta')\kp(\eta_2,\eta'')
d\eta'd\eta''\les\frac{1}{N},
\ees
which yields
\be\label{estimate: I222}
|I_{222}|\leq \frac{1}{N} e^{-\la\al(\tau)}
\sup_{0<s<\tau}\{ e^{\la\al(s)}|w(s)|\}.
\ee

\vskip 0.2 true cm

\textbf{(c) The estimate of $I_{223}$}

\vskip 0.2 true cm

To handle the singular factor  $|\eta_1-\eta'|^{-1}$ in the kernel $\kp(\eta_1,\eta')$,
we choose a smooth function with compact support $k_N$ such that
\bes
\sup_{|p|\leq 3N}\int_{|\eta'|\leq 3N}
|\kp(p,\eta')-k_N(p,\eta')|d\eta' \les \frac{1}{N}
\ees
and then split
\bes
\begin{aligned}
\kp(\eta_1,\eta')\kp(\eta_2,\eta'')
= & \{\kp(\eta_1,\eta')-k_N(\eta_1,\eta')\}\kp(\eta_2,\eta'') \\
& +k_N(\eta_1,\eta')\{\kp(\eta_2,\eta'')-k_N(\eta_2,\eta'')\}\\
& +k_N(\eta_1,\eta')k_N(\eta_2,\eta'').
\end{aligned}
\ees
Using the fact that
\bes
\int_{|\eta''|\leq 3N}|\kp(\eta_2,\eta'')|d\eta''\les 1,
\ees
and
\bes
\int_{|\eta'|\leq 2N}|k_N(\eta_1,\eta')|d\eta'\les 1,
\ees
and noting that the smooth function $k_N$ is bounded by some constant $C_{N}>0$, we obtain that
\be\label{estimate: I223 pre}
|I_{223}|\leq\frac{1}{N} e^{-\la\al(\tau)}
\sup_{0<s<\tau}\{ e^{\la\al(s)}|w(s)|\}
+C_{N} e^{-\la\al(\tau)}J,
\ee
where
\bes
J:=\int_{0}^{\tau}\int_{0}^{\tau_1}\int_{|\eta''|\leq3N}
\int_{|\eta'|\leq2N} e^{\la\al(\tau_2)}
|w(\tau_2,\eta_2,\eta'')|d\tau_2d\tau_1d\eta''d\eta'.
\ees
Choosing
\bes\ba
\de'= & h^2\ka/N^2, \\
\de''= & h^4\ka/N^3,
\ea\ees
in Lemma 2.4 (c). We separate the time intervals $(0,\tau)$
and $(0,\tau_1)$ into the following parts:
\bes
(0,\tau)=\{\cup_k(\tau_{k+1}'+\de',\tau_{k}'-\de')\}\cup
\{\cup_k[\tau_k'-\de',\tau_k'+\de']\},
\ees
\bes
(0,\tau_1)=\{\cup_j(\tau_{j+1}''+\de'',\tau_{j}''-\de'')\}\cup
\{\cup_j[\tau_j''-\de'',\tau_j''+\de'']\}.
\ees
By (d) of Lemma 2.4, we obtain
\bes
J\leq h^2
\sup_{0<s<\tau}\{ e^{\la\al(s)}|w(s)|\}+\sum_{k,j}J_{kj},
\ees
where
\bes
J_{kj}:=\int_{\tau_{k+1}'+\de'}^{\tau_k'-\de'}
\int_{\tau_{j+1}''+\de''}^{\tau_j''-\de''}\int_{|\eta''|\leq3N}
\int_{|\eta'|\leq2N} e^{\la\al(\tau_2)}
|w(\tau_2,\eta_2,\eta'')|d\tau_2d\tau_1d\eta''d\eta'.
\ees
Applying Lemma 2.4 (a) for set $\AehN$, we know that $(y,\eta)\in\AehN$ and $(y_1,\eta_1)\in\AehhN$.
This derives that the time intervals $\Delta\tau'$ and  $\Delta\tau''$ between two adjacent bounces on each trajectory
are
\bes
\Delta\tau'\geq\frac{\ka}{2N^2},\,\,\Delta\tau''\geq\frac{h^2\ka}{2N^3}.
\ees
Moreover, the summations of $k$ and $j$ in (\ref{eq: Y1 and H1}) and (\ref{eq: Y2 and H2}) are finite:
\be\label{estimates: k and j}
k\leq \frac{2N^2\tau}{\ka},\,\,j\leq\frac{2N^3\tau}{h^2\ka}.
\ee
To apply the $L^2$ decay of $u(\tau)$ to derive the bound of $J_{kj}$, we will make a transformation of coordinate: $y_2\mapsto\eta'$, where $y_2=y_2(\tau_2;\tau_1,y_1(\tau_1;\tau,y,\eta),\eta')$. As in Lemma 22 of \cite{MR2679358}, we establish the following result
on $|\det\left(\frac{\p y_2}{\p\eta'}\right)|$.

{\bf Lemma 5.2.} {\it For fixed $k$ and $j$,
$|\det\left(\frac{\p y_2}{\p\eta'}\right)|$ is an analytic function. For any $\ep>0$ small enough, there
are a number $\sigma=\sigma(\ka,N,\ep,k,j)>0$ and an open covering $\cup_{i=1}^mB(\tau_i,y_i,\eta_i;r_i)$ of $(0,\pi/2h)\times\AehN$,
and corresponding open sets $O_{\tau_i,y_i,\eta_i}$ related to $[\tau_{k+1}'+\de',\tau_k'-\de']\times[\tau_{j+1}''+\de'',\tau_{j}''-\de'']
\times\R^3_{\eta'}$ with $|O_{\tau_i,y_i,\eta_i}|<h\ep$, such that
\bes
|\det\left(\frac{\p y_2}{\p\eta'}\right)|\geq
\frac{\sigma}{h^3}>0
\ees
holds for $(\tau,y,\eta)\in(0,\pi/2h)\times\AehN$ and $(\tau_1,\tau_2,\eta')\in O_{\tau_i,y_i,\eta_i}^c\cap
[\tau_{k+1}'+\de',\tau_k'-\de']\times[\tau_{j+1}''+\de'',\tau_{j}''-\de'']
\times\{|\eta'|\leq 2N\}$.
}

{\bf Proof.} We only require to prove Lemma 5.2 for $h=1$. Indeed, for general $h>0$, by choosing new coordinates,
\bes
\theta=h\tau,\, y=y,\, \zeta=\frac{1}{h}\eta,
\ees
then the backward trajectory equations now turn to
\bes\ba
\frac{dy}{d\theta} & =\zeta, \\
\frac{d\zeta}{d\theta} & = -y,
\ea\ees
which just corresponds to the case of $h=1$. Hence, we have the estimates
\bes
|\det\left(\frac{\p y_2}{\p\eta'}\right)|
=\frac{1}{h^3}|\det\left(\frac{\p y_2}{\p\zeta'}\right)|
\geq
\frac{\sigma}{h^3}>0,
\ees
and the open sets $|O_{\tau_i,y_i,\eta_i}|
=h|O_{\theta_i,y_i,\zeta_i}|\leq h\ep$.

Next we assume $h=1$. For fixed $k$, when $\tau_1\in[\tau_{k+1}'+\de',\tau_k'-\de']$,
we have $y_1=y_1(\tau_1;\tau,y,\eta)$ satisfying $(y_1,\eta')\in\AehhN$. So for fixed $j$,
when $\tau_2\in[\tau_{j+1}''+\de'',\tau_j''-\de'']$, $y_2$ is an analytic function of $(\tau_1, \tau_2, \eta')$.
This can be seen from the explicit formula of $y_2$ and Lemma 2.4 (e).

We now start to show that
$|\det\left(\frac{\p y_2}{\p\eta'}\right)|$ is not identically zero. To do so, choose $\tau_2$ such that $|\tau_2-\tau_1|
\leq \frac{h^2\ka}{2N^3}\leq \Delta\tau''$, then
\bes
\begin{aligned}
y_2(\tau_2):= &
Y(\tau_2;\tau_1,y_1,\eta') \\
= &
y_1\cos[h(\tau_2-\tau_1)]
+\frac{\eta'}{h} \sin[h(\tau_2-\tau_1)].
\end{aligned}
\ees
So
\bes
|\det\left(\frac{\p y_2}{\p\eta'}\right)|
=|\frac{\sin[h(\tau_2-\tau_1)]}{h}|^3.
\ees
That is, $y_2$ is an analytic function of $\tau_1,\tau_2,\eta'$ and is not identically zero.
From Lemma 22 of \cite{MR2679358} or p.240 of \cite{MR0257325}, we conclude the
proof of Lemma 5.2.

\vskip 0.2 true cm

Now we derive the estimates of $J_{kj}$. We split the integral in $J_{kj}$
as two parts: including $O_{\tau_i,y_i,\eta_i}$ and including $O_{\tau_i,y_i,\eta_i}^c$ respectively.
Since $|O_{\tau_i,y_i,\eta_i}|<h\ep$, we have
\begin{align}\label{R-2}
& \int_{\tau_{k+1}'+\de'}^{\tau_k'-\de'}
\int_{\tau_{j+1}''+\de''}^{\tau_j''-\de''}\int_{|\eta''|\leq3N}
\int_{|\eta'|\leq2N}1_{O_{\tau_i,y_i,\eta_i}}
 e^{\la\al(\tau_2)}|w(\tau_2,y_2,\eta'')|
d\tau_2d\tau_1d\eta''d\eta' \no\\
&\le h\ep\|w(\tau)\|_{\infty}.
\end{align}
For the second part, we have
\begin{align}\label{R-3}
& \int_{\tau_{k+1}'+\de'}^{\tau_k'-\de'}
\int_{\tau_{j+1}''+\de''}^{\tau_j''-\de''}\int_{|\eta''|\leq3N}
\int_{|\eta'|\leq2N}1_{O_{\tau_i,y_i,\eta_i}^c}
 e^{\la\al(\tau_2)}|w(\tau_2,y_2,\eta'')|
d\tau_2d\tau_1d\eta''d\eta' \no\\
= & \int_{\tau_{k+1}'+\de'}^{\tau_k'-\de'}
\int_{\tau_{j+1}''+\de''}^{\tau_j''-\de''}\int_{|\eta''|\leq3N}
\int_{\O}1_{O_{\tau_i,y_i,\eta_i}^c}
 e^{\la\al(\tau_2)}
\frac{|w(\tau_2,y,\eta'')|}
{|\det(\frac{\p y_2}{\p\eta'})|}
d\tau_2d\tau_1d\eta''dy \no\\
\leq & \frac{h^3C_{\ka,N,\ep}}{\sigma}\int_{0}^{\tau}\int_{0}^{\tau_1}
 e^{\la\al(\tau_2)}\|w(\tau_2)\|_2d\tau_2 d\tau_1 \no\\
\leq & \frac{h^3C_{\ka,N,\ep}}{\sigma}\int_{0}^{\tau}\int_{0}^{\tau_1}
 e^{\la\al(\tau_2)}\|u(\tau_2)\|_2d\tau_2 d\tau_1.
\end{align}
Combining \eqref{R-2}-\eqref{R-3} with the $L^2$ decay of $u$ in (\ref{estimates: u linear}) yields
\be\label{estimate: I223}
|I_{223}|\leq C_{\kappa,N,\epsilon}e^{-\lambda\al(\tau)}.
\ee
Collecting all the estimates (\ref{estimate: I1})-(\ref{estimate: I3}), (\ref{estimate: I21})-(\ref{estimate: I23}),
(\ref{estimate: I221})-(\ref{estimate: I222}) and (\ref{estimate: I223}), we conclude that
\be\label{estimate: w(tau) in A}
\begin{aligned}
\|w(\tau)1_{\AehN}\|_{\infty}\leq & C e^{-\la\al(\tau)}\|w_0\|_{\infty}
+C\left(h^2+\ep+\frac{1}{N}\right)
\sup_{0<s<\tau}\{ e^{\la\al(s)}\|w(s)\|_{\infty}\} \\
& +C\sup_{0<s<\tau}\{ e^{\la\al(s)}\|\nu^{-1}\phi g(s)\|_{\infty}\}
+C_{\ka,N,\ep} e^{-\la\al(\tau)}.
\end{aligned}
\ee

\vskip 0.2 true cm

\textbf{Step 2. Estimate of $\|w(\tau)\|_{\infty}$}

\vskip 0.2 true cm

By expression (\ref{eq: w(tau)}) of $w(\tau)$, for any $(\tau,y,\eta)\in(0,\pi/2h)\times\O\times\R^3$,
$I_1$ and $I_3$ satisfy the same estimates \eqref{estimate: I1} and \eqref{estimate: I3} as in Step 1.
We only need to estimate $I_2$. It follows from direct computation that
\begin{align}\label{R-4}
|I_2|\leq  e^{-\la\al(\tau)}\int_{0}^{\tau}
 e^{\la\al(\tau_1)}|\Kp w(\tau_1)|d\tau_1.
\end{align}
We rewrite $\Kp w(\tau_1)$ as
\bes
\begin{aligned}
\Kp w(\tau_1)   =&\int_{\R^3}\kp(\eta_1,\eta') \{w\cdot(1-1_{\AehN})\}(\tau_1,y_1,\eta') d\eta'\\
& +\int_{\R^3}\kp(\eta_1,\eta') \{w\cdot1_{\AehN}\}(\tau_1,y_1,\eta')d\eta'\\
&:=K_1+K_2.
\end{aligned}
\ees
Similar to the treatments  in Step 1, we have
\begin{align}\label{R-5}
\int_{0}^{\tau}
e^{\la\al(\tau_1)}K_2d\tau_1\le &C e^{-\la\al(\tau)}\|w_0\|_{\infty}
+C\left(h^2+\ep+\frac{1}{N}\right)
\sup_{0<s<\tau}\{ e^{\la\al(s)}\|w(s)\|_{\infty}\}\no\\
& +C\sup_{0<s<\tau}\{ e^{\la\al(s)}\|\nu^{-1}\phi g(s)\|_{\infty}\}
+C_{\ka,N,\ep} e^{-\la\al(\tau)}.
\end{align}
To estimate $\int_{0}^{\tau}
e^{\la\al(\tau_1)}K_1d\tau_1$, we need to consider the following three cases of the integration
variable $\eta'$: $|\eta'|<4h$, $|\eta'|>2N$ and  $4h\leq|\eta'|\leq 2N$.
In this case,
\bes
\ba
K_1= & \left(\int_{|\eta'|<4h}+\int_{|\eta'|>2N}+\int_{4h\leq|\eta'|\leq2N}\right)
  \kp(\eta_1,\eta') \{w\cdot(1-1_{\AehN})\}(\tau_1,y_1,\eta') d\eta'  \\
  := & K_{11}+K_{12}+K_{13}.
\ea
\ees
By the measure $|\{\eta'\in\R^3:|\eta'|<4h\}|\lesssim h^3$, we have
\be\label{estimates: K11}
  |K_{11}|\lesssim h^3.
\ee
In addition, it follows from Lemma 3.2 and (\ref{estimate: k}) that for  large $N$,
\be\label{estimates: K12}
  |K_{12}|\lesssim 1/N.
\ee
Note that by the definition of $\AehN$ in (\ref{def: set A}), one has
\bes\ba
  & \{\eta'\in\R^3: 4h\leq|\eta'|\leq2N\}\cap\{\eta'\in\R^3: (y_1,\eta')\in\R^3-\AehN\} \\
 = & \{\eta'\in\R^3:4h\leq|\eta'|\leq2N \text{ and } |\eta'|^2+h^2|y_1|^2-|y_1\times\eta'|^2<h^2+\kappa^2\} .
\ea\ees
This means that
\bes
\begin{aligned}
	\ka^2> & |\eta'|^2+h^2|y_1|^2-|y_1\times\eta'|^2-h^2 \\
	= & (|\eta'|^2-h^2)(1-|y_1|^2)+|\eta'\cdot y_1|^2,
\end{aligned}
\ees
which yields $|y_1|^2\geq 1-\ka^2/h^2$ and $|\eta'\cdot y_1|\leq\ka$. For fixed $y_1\in\Omega$, we arrive at
\bes
  |\{\eta'\in\R^3: 4h\leq|\eta'|\leq2N \text{ and } |\eta'|^2+h^2|y_1|^2-|y_1\times\eta'|^2<h^2+\kappa^2\}|\lesssim N^2\kappa.
\ees
Choosing $\kappa<1/N^3$, then we have
\be\label{estimates: K13}
  |K_{13}|\lesssim 1/N.
\ee

Collecting (\ref{estimates: K11})-(\ref{estimates: K13}) and \eqref{R-5},
and Choosing $N$ large enough and $\ep>0$ small enough, we eventually get
\bes
\sup_{0<s<\tau}\{ e^{-\la\al(s)}\|w(\tau)\|_{\infty}\}
\lesssim\|w_0\|_{\infty}
+\sup_{0<s<\tau}
\{ e^{\la\al(s)}\|\nu^{-1}\phi g(s)\|_{\infty}\}.
\ees
Then we complete the proof of Proposition 5.1. \qquad\qquad \qquad \qquad \qquad \qquad \qquad \qquad
\qquad \qquad  $\square$

\section{Proofs of Theorem 1.1 and 1.2}
Set $w=\phi u$. Then we have
\begin{equation}\label{eq: w nonlinear}
\Lat w+\mut(y)\cos^2(h\tau)\nu(\eta)w \\
=\mut(y)\cos^2(h\tau)\Kp w+\mut^{1/2}(y)\cos^2(h\tau)\Gp(w,w).
\end{equation}
We will prove Theorem 1.1 by the standard Picard iteration.

\vskip 0.2 true cm

{\bf Proof of Theorem 1.1.} The proof will be divided into the following four steps.

\vskip 0.2 true cm

{\bf Step 1. Existence of solution $w$ to \eqref{eq: w nonlinear}}

\vskip 0.2 true cm

Let $w^0\equiv 0$. For $m\geq 0$, we define the following iteration:
\begin{equation}\label{eq: iteration}
\{\Lat+\mut(y) \cos^2(h\tau)\nu(\eta)\}\wmp \\
=\mut(y) \cos^2(h\tau)\Kp\wm+\mut^{1/2}(y) \cos^2(h\tau)\Gp(\wm,\wm),
\end{equation}
with the initial-boundary data
\begin{equation}\label{data for wmp}
\begin{cases}
\wmp(0,y,\eta)=w_0(y,\eta),&
\text{for $y\in \O$, $\eta\in\R^3$},\\
\wmp(\tau,y,\eta)=\wmp(\tau,y,\eta-2(\eta\cdot n)n),&
\text{for $\eta\cdot n<0$, $y\in\p\O$}.
\end{cases}
\end{equation}
As in Section 5,  a mild solution $w(\tau,y,\eta)$
of (\ref{eq: iteration}) with (\ref{data for wmp}) can be explicitly constructed.
Moreover, by Proposition 5.1 and Lemma 3.4 for $\Gp$, we have
\begin{align}\label{R-6}
\|\wmp(\tau)\|_{\infty}\leq C e^{-\la\al(\tau)}\|w_0\|_{\infty}
+C e^{-\la\al(\tau)}
\sup_{0<s<\tau}\{ e^{\la\al(s)}\|\wm(s)\|_{\infty}^2\}.
\end{align}
Assume that $\|\wm(\tau)\|_{\infty}\leq 2C e^{-\la\al(\tau)}\|w_0\|_{\infty}$, which is true for
$m=0$, then it follows from \eqref{R-6} that
\begin{align}\label{R-7}
\|\wmp(\tau)\|_{\infty}\leq C e^{-\la\al(\tau)}\|w_0\|_{\infty}
+(4C^2\|w_0\|_{\infty})C e^{-\la\al(\tau)}\|w_0\|_{\infty}.
\end{align}
If the initial data satisfies that $\|w_0\|_{\infty}\leq 1/4C^2$, we then conclude that
\bes
\|\wmp(\tau)\|_{\infty}\leq 2C e^{-\la\al(\tau)}\|w_0\|_{\infty}.
\ees
Thus, by induction method, we have that for all $m\geq 0$,
\be\label{estimate: wm}
\|\wm(\tau)\|_{\infty}\leq 2C e^{-\la\al(\tau)}\|w_0\|_{\infty}.
\ee
On the other hand, setting $\gmp=\wmp-\wm$ yields
\be\label{eq: iteration of diffence}
\ba
\{\Lat+\mut(y)\cos^2(h\tau)\nu(\eta)\}\gmp  = & \mut(y)\cos^2(h\tau)\Kp\gm +\mut^{1/2}(y)\cos^2(h\tau)\{\Gp(\wm,\wm)
\\
& -\Gp(\wmm,\wmm)\}.
\ea\ee
Note that
\begin{align*}
&|\Gamma_{\phi}(w^{m},w^{m})-\Gamma_{\phi}(w^{m-1},w^{m-1})|\\
=&|\Gamma_{\phi}(w^{m}-w^{m-1},w^{m})+\Gamma_{\phi}(w^{m-1},w^m-w^{m-1})|\\
\leq & C\nu(\eta)\|w^m-w^{m-1}\|_{\infty}
\left(\|w^m\|_{\infty}+\|w^{m-1}\|_{\infty}\right).\\
\leq & C\nu(\eta)\|w_0\|_{\infty}\|g^m\|_{\infty}.
\end{align*}
Then by the analogous estimate in Proposition 5.1 for (\ref{eq: iteration of diffence}), one has
\bes
\|\gmp(\tau)\|_{\infty}\leq  C e^{-\la\al(\tau)}\|w_0\|_{\infty}
\sup_{0<s<\tau}\{ e^{\la\al(s)}\|\gm(s)\|_{\infty}\}.
\ees
By assuming $\|w_0\|_{\infty}\leq 1/4C^2$ as in the above, we arrive at
\bes
\|\gmp(\tau)\|_{\infty}\leq \frac{1}{4C}\|\gm(\tau)\|_{\infty}.
\ees
Thus there exists a function $w(\tau)$ such that $\wm(\tau)\to w(\tau)$ in $L^{\infty}$
and $w$ is a mild solution to \eqref{eq: w nonlinear} with \eqref{data for wmp}.

\vskip 0.2 true cm

{\bf Step 2. Uniqueness  of solution $w$ to \eqref{eq: w nonlinear}}

\vskip 0.2 true cm

Assume that there is another solution $\bw$ to
\eqref{eq: w nonlinear} with the same initial-boundary data as $w$, and also assume that
$\sup_{\tau}\{ e^{\la\al(\tau)}\|\bw(\tau)\|_{\infty}\}$ is small. Then
\bes\ba
\{\Lat+\mut(y)\cos^2(h\tau)\nu(\eta)\}\{w-\bw\} = & \mut(y)\cos^2(h\tau)\Kp\{w-\bw\} \\
  & +\mut^{1/2}(y)\cos^2(h\tau)\{\Gp(w,w)-\Gp(\bw,\bw)\}
\ea\ees
with the vanishing initial data for $w-\bw$.  As in Step 1,
we can derive that $\|\{w-\bw\}(\tau)\|_{\infty}\equiv 0$. Therefore,
the uniqueness  of solution $w$ to \eqref{eq: w nonlinear} is shown.

\vskip 0.2 true cm

{\bf Step 3. Positivity  of solution $f$ to \eqref{Boltzmann eq}}

\vskip 0.2 true cm

Let $f^0=f_0$, by \eqref{Boltzmann eq}, we solve $f^{m+1}$ for $m\ge0$ as follows
\begin{equation}\label{eq: iteration of f}
\Lat f^{m+1}+\cos^2(h\tau)\nu(f^m)f^{m+1}
=\cos^2(h\tau) Q_1(f^m,f^m),
\end{equation}
where
\begin{equation*}
\nu(f^m)=\int_{\R^3\times S^2}|(\eta_*-\eta)\cdot\omega|f^m(\eta_*)d\eta_*d\omega
\end{equation*}
and
\begin{align*}
Q_1(f^m,f^m)=& Q(f^m,f^m)-2\nu(f^m)f^m\\
=& \int_{\R^3\times S^2}
|(\eta_*-\eta)\cdot\omega|f^m(\eta')f^m(\eta'_*)d\eta_*d\omega.
\end{align*}
Set $f^m=M+M^{1/2}u^m$ and $f^{m+1}=M+M^{1/2}u^{m+1}$. Then it follows from (\ref{eq: iteration of f}) that
\bes\ba
\{\Lat+\mut(y)\cos^2(h\tau)\nu(\eta)\}u^{m+1} = & \mut(y)\cos^2(h\tau)Ku^m \\
 & +\mut^{1/2}(y)\cos^2(h\tau)\{\Gamma_1(u^m,u^m)-\Gamma_3(u^m,u^{m+1})\}.
\ea\ees
One can check that $w^m=\phi u^m$  converges in $L^{\infty}$ to the solution $w$
of \eqref{eq: w nonlinear} as in Step 1.

Assume $f^m\geq0$. Let
\begin{equation*}
\alpha(\tau)=\cos^2(h\tau)\nu(f^m).
\end{equation*}
Then by integrating along the backward trajectories $(Y(\tau), H(\tau))$ of \eqref{eq: iteration of f}, we have
that for $(\tau,y,\eta)\in S_0\times\R^3$,
\bes\ba
f^{m+1}(\tau,y,\eta)= & f_0(Y(0),H(0))e^{-\int_0^{\tau}\alpha(s)ds} \\
& +\int_0^{\tau}\cos^2(h\tau) Q_1(f^m,f^m)(s,Y(s),H(s)) e^{-\int_t^{\tau}\alpha(s)ds}dt.
\ea\ees
This derives $f^{m+1}\geq0$ and we further deduce that the solution $f\geq0$ of \eqref{Boltzmann eq}
by the uniqueness in Step 2.

\vskip 0.2 true cm

{\bf Step 4. Continuity  of solution $f$ to \eqref{Boltzmann eq}}

\vskip 0.2 true cm

The continuity of the solution $f$ is obvious since we have obtained the continuity of the backward trajectory in Lemma 2.5.
\qquad \qquad \qquad \qquad \qquad \qquad \qquad \qquad \qquad \qquad \qquad \qquad \qquad \qquad
\qquad $\square$

\vskip 0.3 true cm

Based on Theorem 1.1, we start to prove Theorem 1.2.

\vskip 0.2 true cm

{\bf Proof of Theorem 1.2.}
The initial data $f_0$ can be reformulated as
\begin{equation*}
f_0=M+M^{1/2}u_0,
\end{equation*}
where
\begin{equation*}
u_0=\left( e^{h^2|y|^2/2}- e^{-h^2|y|^2/2}\right) e^{-|\eta|^2/2}+ e^{-h^2|y|^2/2}\tilde{u_0},
\end{equation*}
and
\begin{equation*}
|u_0|\leq 2h^2 e^{-|\eta|^2/2}+|\tilde{u_0}|.
\end{equation*}
When $\|\phi\tilde{u_0}\|_{\infty}<\epsilon$ and $h<\epsilon^{1/2}$, we have
\begin{equation*}
\|\phi u_0\|_{\infty}\leq C\epsilon.
\end{equation*}
In this case, all the assumptions in Theorem 1.1 are fulfilled. Therefore, by Theorem 1.1,
there exists a unique mild solution $f=M+M^{1/2}u\geq0$ to problem (\ref{eq: f}) with (\ref{data: f}).
Going back to the original coordinates $(t,x,\xi)$, the perturbation solution $u$ satisfies
\begin{equation*}
\|\phi u(t)\|_{\infty}\leq C\|\phi u_0\|_{\infty}\leq C\epsilon.
\end{equation*}
Next we derive the decay property of the mass density $\rho$ of gases. Note that
\begin{align*}
\rho_1(t,x):= & \int_{\R^3}M(t,x,\xi)d\xi\\
= & \frac{\pi^{3/2}}{R(t)^3}\exp{\left(-\frac{h^2|x|^2}{R(t)^2}\right)},\\
\rho_2(t,x):=&\int_{\R^3}|M(t,x,\xi)^{1/2}u(t,x,\xi)|d\xi\\
\leq & \|\phi u\|_{\infty}\int_{\R^3}M(t,x,\xi)^{1/2}d\xi\\
= & \|\phi u\|_{\infty}\frac{(2\pi)^{3/2}}{R(t)^3}\exp{\left(-\frac{h^2|x|^2}{2R(t)^2}\right)}.
\end{align*}
Hence, when $\epsilon>0$ is small, by $\rho(t,x)=\rho_1(t,x)+\rho_2(t,x)$ we have
\begin{equation*}
\frac{1}{C_0R(t)^3}\leq\rho(t,x)\leq\frac{C_0}{R(t)^3},
\end{equation*}
where $C_0>1$ is a constant. With respect to the continuity of the solution $f$ to problem (\ref{eq: f}) with (\ref{data: f}),
we only need to change $\ga_{01}$ in Theorem 1.1 to the corresponding set in coordinates $(t,x,\xi)$. Consequently,
we complete the proof of Theorem 1.2. \qquad \qquad \qquad \qquad \qquad \qquad \qquad \qquad \qquad $\square$

\appendix

\section{The study on elliptic system (\ref{eq: elliptic partial tau Phib}) and (\ref{eq: vphi b})}
In this section, we will prove the existence of the solution to the elliptic system (\ref{eq: elliptic partial tau Phib}) and (\ref{eq: vphi b}).
For a vector function $b(y)$ defined in $\O=\{|y|<1\}$, consider the following elliptic system of $\vp=
(\vp_1, \vp_2, \vp_3)$:
\begin{equation}\label{eq: vphi bb}
\begin{cases}
-\Lap  \vphi=   b\qquad   & \text{ in $\O$}, \\
\vphi\cdot n=   0 \qquad & \text{ on $\p\O$}, \\
\p_n\vphi=  (\p_n\vphi\cdot n)n  \qquad & \text{ on $\p\O$}.
\end{cases}
\end{equation}
The second boundary condition in \eqref{eq: vphi bb}  can be rewritten as
\begin{align}\label{A-1}
\left(
\begin{array}{ccc}
1-(n^1)^2 & -n^1n^2 & -n^2n^3 \\
-n^2n^1 & 1-(n^2)^2 & -n^2n^3 \\
-n^3n^1 & n^3n^2 & 1-(n^3)^2 \\
\end{array}
\right)
\left(
\begin{array}{c}
\p_n\vphi^1 \\
\p_n\vphi^2 \\
\p_n\vphi^3 \\
\end{array}
\right)
=0.
\end{align}
Note that  the rank of the coefficient matrix in \eqref{A-1} is $2$.

Define the following  Banach space $\mathcal{V}$
\bes
\msV:=\{\psi\in(H^1(\O))^3:\,\psi\cdot n=0 \text{ on $\p\O$}\}
\ees
with the norm
\bes
\|\psi\|_{\msV}^2=\|\psi\|_{H^1}^2=\|\psi\|^2+\|\nabla\psi\|^2.
\ees
Obviously,
\bes
(H^1_0(\O))^3\subset\msV\subset (H^1(\O))^3,
\ees
and the dual space of $\msV^*$ satisfies
\bes
(H^{-1}_0(\O))^3\subset\msV^*\subset (H^{-1}(\O))^3.
\ees
For $\vphi, \psi\in\msV$ and $\vphi\in H^2$, we have
\be\label{eq: ellitic Green}
-(\Delta\vphi,\psi)=(\nabla\vphi,\nabla\psi)
-\int_{\p\O}\psi\cdot\p_n\vphi.
\ee
By the first boundary condition of $\psi\in\msV$ in (\ref{eq: vphi bb}), the boundary term
in (\ref{eq: ellitic Green}) turns to
\bes\ba
\int_{\p\O}\psi\cdot\p_n\vphi= &
\int_{\p\O}[\psi-(\psi\cdot n)n]\cdot\p_n\vphi \\
= & \int_{\p\O}\psi\cdot\p_n\vphi
-\int_{\p\O}(\psi\cdot n)(n\cdot\p_n\vphi) \\
= & \int_{\p\O}\psi\cdot[\p_n\vphi-(n\cdot\p_n\vphi)n].
\ea\ees
This means that the boundary term  in (\ref{eq: vphi bb}) vanishes if $\vphi$ satisfies the second
boundary condition in (\ref{eq: vphi bb}).
Therefore, we can define the operator $\msL:\,\msV\to\msV^*$ as follows
\bes
(\msL\vphi,\psi)=(\nabla\vphi,\nabla\psi).
\ees
It is easy to check that $\msL$ is a bounded self-adjoint operator. We now show that $\msL$ has a closed range.
In fact, since
\bes\ba
\|\psi\|_{\msV}^2= & \|\psi\|^2+\|\nabla\psi\|^2 \\
= & (\msL\psi,\psi)+\|\psi\|^2 \\
\leq & \|\msL\psi\|_{\msV^*}\|\psi\|_{\msV}+\|\psi\|^2,
\ea\ees
we have
\bes
\|\psi\|_{\msV}^2\leq C(\|\msL\psi\|_{\msV^*}^2+\|\psi\|^2).
\ees
Together with the fact that the mapping $\mathrm{id}:\,\msV\to L^2$ is compact, we know that $\msL$
has a closed range by Proposition 6.7 in Appendix A of \cite{MR2744150}.

Next we  prove that $\msL$  is a one-to-one and onto mapping. For $\vphi_0\in\text{Ker}\msL$, due to
\bes
0=(\msL\vphi_0,\vphi_0)=(\nabla\vphi_0,\nabla\vphi_0)
=\|\nabla\vphi_0\|^2,
\ees
$\vphi_0$ is a constant vector in $\O$. This, together with the boundary condition $\vphi_0\cdot n=0$,
yields $\vphi_0=0$. So $\msL$ is a one-to-one mapping. In addition, since $\msL$ is self-adjoint and has a
closed range, we obtain the range $\mathcal{R}(\msL)=(\text{Ker}\msL)^{\perp}=\msV^*$. This means
that $\msL$ is
also an onto mapping. Consequently, the bounded bijective linear operator $\msL$ has a bounded inverse $\msL^{-1}$.
So problem \eqref{eq: vphi bb} is uniquely solved by $\vphi=-\msL^{-1}b\in\msV$ in a weak sense.
On the other hand, the existence of a classical solution $\vp$  to problem \eqref{eq: vphi bb} can
be obtained by the standard methods in Chapter 5 of \cite{MR2744150}. We omit the proof here.

\vskip 0.4 true cm

{\bf Acknowledgement.} {\it The authors  wish to express their gratitude to Professor Yang Tong,
the City University of Hong Kong, for his
interests in this problem and many very fruitful discussions in the past. In particular,
Professor Yang Tong gave many suggestions and comments in this topic.}

\vskip 0.3 true cm

\bibliographystyle{plain}
%\bibliography{reference}

\end{document}